\documentclass[10pt]{article}
\usepackage{graphicx}
\usepackage{amsfonts}
\usepackage{bbm}
\usepackage{mathrsfs}
 \usepackage{mathrsfs,amsmath,amssymb}
\usepackage[dvips]{color}
 \usepackage{amssymb}
 \usepackage[dvips]{color}
\allowdisplaybreaks
 \renewcommand{\baselinestretch}{1}
 
 \newtheorem{remark}{\noindent\mbox{Remark}}[section]

 \newtheorem{theorem}{\noindent\mbox{Theorem}}[section]
 \newtheorem{proposition}{\noindent\mbox{Proposition}}[section]
 \newtheorem{corollary}{\noindent\mbox{Corollary}}[section]
 \newtheorem{definition}{\noindent\mbox{Definition}}[section]

 \def\bq{\begin{equation}}
 \def\eq{\end{equation}}
 \def\eqn{\end{eqnarray}}
 \def\bqn{\begin{eqnarray}}
 
 \def\qed{\hfill$\Box$\medskip}
 \def\rto{\rightarrow\infty}
 \def\z{\left}
 \def\y{\right}
 \def\no{\nonumber}

 \def\mcp{\mathcal{P}}
 \def\om{\overline{\omega}}
 
\begin{document}
\noindent{\Large\bf { Intrinsic} branching structure {
within} random walk on $\mathbb{Z}$\footnote{
 The project is partially supported by National
Nature Science Foundation of China (Grant No.10721091) and NCET(
No.\,05-0143).}} 

\noindent{
  Wenming Hong$^{a}$ \&
Huaming Wang$^{b}$
}

\vspace{0.2 true cm}
\renewcommand{\baselinestretch}{1.3}\baselineskip 12pt
\noindent{\footnotesize\rm  $^a$ School of Mathematical Sciences \&  Key
Laboratory of Mathematics and Complex Systems, Beijing Normal
University, Beijing 100875, P.R.China;\\%
 $^b$ Business College, Beijing Union University, Beijing 100025, P.R.China \\
 (email: wmhong@bnu.edu.cn, huamingking@mail.bnu.edu.cn)
\vspace{4mm}}%

\begin{center}
\begin{minipage}[c]{15cm}
\begin{center}\textbf{Abstract}\end{center}
In this paper, we reveal the branching structure for a
non-homogeneous random walk with bounded jumps. The ladder time
$T_1,$ the first hitting time of $[1,\infty)$ by the walk starting
from $0,$ could be expressed in terms of a non-homogeneous multitype
branching process. As an application of the branching structure,
{ we prove a law of large numbers of random walk in random
environment with bounded jumps and  specify the explicit invariant
density for the Markov chain of ``the environment viewed from the
particle" .The invariant density and the limit velocity could be
expressed explicitly   in terms of the environment.} \vspace{0.2cm}

\mbox{}\textbf{Keywords:}\quad random walk, branching process, random environment, invariant density.\\
\mbox{}\textbf{Mathematics Subject Classification}:  Primary 60J80;
secondary 60G50.
\end{minipage}
\end{center}

\section{Introduction}
\subsection{Background}
We study random walk $\{X_n\}$ with bounded jumps in this paper.
More precisely, let $\Lambda=\{-L,...,R\}/\{0\}.$ At position $x$
the walk jumps with probability $\omega_x(l)$ to $x+l$ for $l\in
\Lambda.$ Of course, $\sum_{l\in\Lambda}\omega_x(l)=1$ and
$\omega_{x}(l)\ge 0,\ l\in \Lambda. $ We call
$\{X_n\}$ the (L-R) random walk hereafter.

Let $T_1$ be the time the walk hits the positive half line
$(0,\infty)$ for the first time. { As we know, it is a
fundamental task to characterize $T_1$,  which} plays an important
roles in studying the  (L-R) random walk, for example, the
recurrence vs transience, law of large numbers, central limits
theorem and large deviations principle et al.

Indeed, for $T_1,$ even some simple questions are difficult to answer. For
example, what is the distribution of $T_1?$ Or further, what are the
moments of $T_1.$

For some simple setting, these questions are known but are open for
some more general setting.

For simple random walk, that is, $L=R=1$ and
$\omega_x(1)=1-\omega_{x}(1)=p\in[0,1]$ for all $x\in\mathbb{Z}.$
The distribution of $T_1$ could be find by Reflection Principle. We
refer the reader to Strook \cite{s05}(2005) for the specific
calculation.

For (1-1) random walk  with non-homogeneous transition
probabilities, that is, $\omega_x$ depends on $x,$ it has been
revealed in Kesten-Kozlov-Spitzer \cite{kks75}(1975) that $T_1$
could be expressed in term of a non-homogeneous Galton-Watson
branching process. This fact enables them to prove a nice stable limit
theorem for nearest neighbor Random Walk in Random Environment (RWRE
hereafter). If $\max\{L,R\}>1,$ that is, for the non-nearest
neighbor random walk, things get very different.

For (L-1) random walk, Hong-Wang \cite{hw10}(2010) revealed its
branching structure, that is, $T_1$ could be express by a multitype
($L$-type) branching process. Using the branching structure,
Hong-Wang also proved a stable limit theorem, partially generalizing
the results, for (1-1) RWRE, of Kesten-Kozlov-Spitzer \cite{kks75}(1975)  to (L-1) RWRE which is
supposed to be transient to the right.

One may ask naturally  what the branching structure for (1-R) random
walk is. This question is answered in Hong-Zhang \cite{hzh10}(2010).
The authors decomposed the random walk path to revealed that $T_1$
for (1-R) random walk could be expressed by a non-homogeneous
$(1+2+...+R)$-type branching process and by this fact they also
proved a law of large numbers for (1-R) RWRE by a method known as
``the environment viewed from particles".

One should note that, both the above mentioned  \cite{hw10} and
\cite{hzh10} treat the case $\min\{L,R\}=1.$  The walk is requested
to be nearest neighbor at least in one side. The main purpose of
this paper { is to consider the general situation:
$\min\{L,R\}\ge 1.$ We restrict ourselves  to $L=R=2$ to explain the
the idea. We reveal the intrinsic  branching structure within the
(L-R) random walk and give some applications.}

Next we define the precise model of (L-R) random walk.
\subsection{The model}
Fix $L,R\ge1.$ Set $\Lambda=\{-L,...,R\}/\{0\}.$ For $i\in\mathbb{Z},$ let $\omega_{i}=(\omega_i(l))_{l\in\Lambda}$ be a probability measure on $i+\Lambda,$ that is $\sum_{l\in \Lambda}\omega_i(l)=1,$ and $\omega_{i}(l)\ge 0$ for all $l\in\Lambda.$ Set $\omega=\{\omega_i\}_{i\in \mathbb{Z}},$ which will serve as the transition probabilities of the random walk. Let $\{X_n\}_{n\ge0}$ be a Markov chain with initial value $X_0=x$ and transition probabilities
$$P_\omega(X_{n+1}=i+j\big |X_n=i)=\omega_i(j),\ j\in\Lambda.$$
We call $\{X_n\}$ the {\it(L-R) random walk with non-homogeneous transition probabilities}.
Throughout this paper, we use $P_\omega^x$ to denote the law induced by the random walk $\{X_n\}.$

In the remainder of the paper, in order to avoid the heavy notations, we consider the case $L=2,R=2,$ that is, the (2-2) random walk. The idea for treating the general (L-R) random walk is basically the same as (2-2) setting. Also except otherwise stated, we always assume the random walk starts from $0.$

For the above defined (2-2) random walk $\{X_n\},$ set $T_0=0$ and define recursively
\begin{equation}\no
  T_k=\inf[n\ge0:X_n>X_{T_{k-1}}]
\end{equation} for $k\ge1.$
We call the stopping times $T_k,$ $k\ge1,$ the {\it ladder times} of the random walk.

 Especially one sees by the definition that  $$T_1=\inf[n\ge0:X_n> 0]$$ is the hitting time of $[1,\infty)$ by the walk,  and $X_{T_1}$ is also a random variable with two possible values $X_{T_1}=1$ or $X_{T_1}=2.$ This is the reason why we call $T_k,$ $k\ge1,$ the ladder times. The distribution of $X_{T_1},$ the exit probabilities of the walk from $(-\infty,0],$ is also a question we concern below.

 The main purpose of this paper is to count exactly how many steps the walk spends to exit successfully from $(-\infty,0]$ (the total steps of the walk before $T_1.$) and give some applications. Next we state the main results.
 \subsection{The main results}
In order to count exactly the steps of the walk  before $T_1,$ we define three types of excursions.
\begin{definition}\mbox{}

 \textbf{a)} We call excursions of the form $\{X_k=i,X_{k+1}=i-1,X_{k+2}\le
i-1,...,X_{k+l}\le i-1,X_{k+l+1}\ge i\}$ type-$\mathcal{A}$
excursions at $i.$ Corresponding to the three kinds of possible last
step of type-$\mathcal{A}$ excursions at $i,$ say, $\{i-1\rightarrow i\},$  $\{i-2\rightarrow i\}$ and $\{i-1\rightarrow i+1\},$
 we classify type-$\mathcal{A}$ excursions at $i$ into three
sub-types $\mathcal A_{i,1},$ $\mathcal A_{i,2}$ and $\mathcal A_{i,3}.$

\textbf{b)} We call excursions of the form $\{X_k=i,X_{k+1}=i-2,X_{k+2}\le
i-1,...,X_{k+l}\le i-1,X_{k+l+1}\ge i\}$ type-$\mathcal{B}$
excursions at $i.$ Corresponding to the three kinds of possible last
step of type-$\mathcal{B}$ excursions at $i,$ say, $\{i-1\rightarrow i\},$  $\{i-2\rightarrow i\}$ and
$\{i-1\rightarrow i+1\},$ we classify type-$\mathcal{B}$ excursions at $i$ into three
sub-types $\mathcal B_{i,1},$ $\mathcal B_{i,2}$ and $\mathcal B_{i,3}.$

\textbf{c)} We call excursions of the form $\{X_k=i+1,X_{k+1}=i-1,X_{k+2}\le
i-1,...,X_{k+l}\le i-1,X_{k+l+1}\ge i\}$ type-$\mathcal{C}$
excursions at $i.$ Corresponding to the three kinds of possible last
step of type-$\mathcal{C}$ excursions at $i,$ say, $\{i-1\rightarrow i\},$  $\{i-2\rightarrow i\}$ and
 $\{i-1\rightarrow i+1\},$ we classify type-$\mathcal{C}$ excursions at $i$ into three
sub-types $\mathcal C_{i,1},$ $\mathcal C_{i,2}$ and $\mathcal C_{i,3}.$
\end{definition}
\begin{center}
\includegraphics[height=6cm]{excursiontype.1}
\end{center}

Define $$A_{i,j}=\#\{\mathcal A_{i,j} \text{ excursions before } T_1\},$$
$$B_{i,j}=\#\{\mathcal B_{i,j} \text{ excursions before } T_1\},$$
$$C_{i,j}=\#\{\mathcal C_{i,j} \text{ excursions before } T_1\},$$ for $i\le 0$ and $j=1,2,3,$
where $``\#\{\}"$ means the number of the elements in some set.

We aim at counting exactly all steps by the walk before $T_1.$ For this purpose, define
\begin{equation}\label{ui0}
  U_{i}=(A_{i,1},A_{i,2},A_{i,3},B_{i,1},B_{i,2},B_{i,3},C_{i,1},C_{i,2},C_{i,3}),
\end{equation} being the the total numbers of different excursions at $i$ before time $T_1.$
Then we have the following fact.
\begin{theorem}\label{tu}
  Suppose that $\limsup_{n\rto}X_n=\infty.$ Then
  \begin{equation}\label{tue}T_1=1+\sum_{i\le0}U_i(2,2,1,1,1,0,2,2,1)^T,\end{equation} where for vector $\mathbf{v}\in\mathbb{R}^9,$ $\mathbf{v}^T$ denotes the transposition of ${\bf v}.$
\end{theorem}The proof of Theorem \ref{tu} will be given in Section \ref{prtu} below.
\begin{remark}
  Because the (2-2) random walk we considering is non-nearest neighbor, one can not give the exact distribution of $T_1$ in general.
But we find in this paper that the process $\{U_i\}_{i\le 0}$ defined in  (\ref{ui0}) is a non-homogeneous multitype branching process. This fact together with (\ref{tue}) enables us to study $T_1$ by the properties of branching processes.
\end{remark}
We attach first  an ancestor to the branching process. The walk starts from
$0.$ But before $T_1,$ there is no jump down from above $1$ to $0$
by the walk. One can imagine that there is a step by the walk from
$1$ to $0$ before it starts from $0$ (One can also imagine that this
step is from $2$ to $0.$ But this makes no difference.), that is,
set $X_{-1}=1.$ Adding this imaginary step, the path $\{X_{-1}=1,
X_0=0,X_1,...,X_{T_1}\}$ forms a type-$\mathcal{A}$ excursion at $1$  such that with probability $1,$
$$A_{1,1}+A_{1,2}+A_{1,3}=1.$$
Then one defines $U_1$ as in (\ref{ui0}). But since there is no  $\mathcal{B}_{1,j}$ and $\mathcal C_{1,j},$ $j=1,2,3,$ excursions, $U_1$ has only three possible values, that is, $U_1=  \mathbf e_1,$ $U_1=  \mathbf e_2$ or $U_1= \mathbf e_3.$

We can treat $U_1$ as some particle immigrates in the system  and call it ``immigration" throughout.
\begin{theorem}\label{thbranc}Suppose that $\limsup_{n\rto}X_n=\infty.$ Then
   $\{U_{i}\}_{i\le 1}$ is a 9-type
non-homogeneous branching processes with immigration distribution as in (\ref{inia}), (\ref{inib}) and (\ref{inic}) below,  and offsprings distributions as in  (\ref{oaa}-\ref{oad}), (\ref{oba}-\ref{obd}), (\ref{oca}-\ref{ocd}) and (\ref{oda}-\ref{ode}) below.
\end{theorem}
The proof of Theorem \ref{thbranc} will be given in the long Section \ref{bran} below.
\begin{remark}
  \begin{itemize} \item[\textbf{1)}] While $\omega_{i}(-2)=0$ for all $i\in\mathbb Z,$ Theorem \ref{thbranc} reveals the branching structure of (1-2) random walk and while  $\omega_{i}(2)=0$ for all $i\in\mathbb Z,$ Theorem \ref{thbranc} reveals the branching structure of (2-1) random walk. Therefore the branching structure in Theorem \ref{thbranc} contains both the branching structures of Hong-Wang \cite{hw10} and Hong-Zhang \cite{hzh10}. For the details, see Remark \ref{deg} below.
  \item[\textbf{2)}] The authors found that, one could simplify the branching structure. In fact, a $6$-type branching process is enough to count exactly all steps by the walk before $T_1.$ However we still use a $9$-type branching process, since it is more understandable and each of the 9 types of particles corresponds to specific jump of the walk. For details, see also Remark \ref{deg} below.
      \end{itemize}
\end{remark}
We give an example to test the branching structure (the above Theorem \ref{tu} and Theorem \ref{thbranc}) in Section \ref{test}.
We consider the degenerated $\omega,$ that is, for all $i\in\mathbb{Z},$ $\omega_i(1)=p_1,$ $\omega_i(2)=p_2,$ $\omega_i(-1)=q_1$ and $\omega_{i}(-2)=q_2$ with $p_1+p_2+q_1+q_2=1,$ $p_1,p_2,q_1\ge 0$ and $q_2>0.$
In this case, let
  \begin{equation}\label{mc}\no
  M=\left(
    \begin{array}{ccc}
      -\frac{q_1+q_2}{q_2} & \frac{p_1+p_2}{q_2} & \frac{p_2}{q_2} \\
     1  &  0& 0\\
      0& 1 & 0 \\
    \end{array}
  \right).
\end{equation}
Suppose that $E_\omega^0(X_1)=p_1+2p_2-q_1-2q_2>0$ implying that $\lim_{n\rto}X_n=\infty$ and let $f,$ $g,$ $h$ be the three eigenvalues of $M$ such that $|f|>|g|>|h|.$ Then we show that $$P_\omega^0(X_{T_1}=1)=1+h \text{ and }P_\omega^0(X_{T_1}=2)=-h, $$ where $h<0$ follows from $p_1+2p_2-q_1-2q_2>0$ (See Section \ref{test} for details.).
Then one has that \begin{equation}\label{ext1}
  E_\omega^0(X_{T_1})=1-h.
\end{equation}

On the other hand, as $\omega$ is degenerated, one follows from  Theorem \ref{tu} and Theorem \ref{thbranc} that
\begin{equation}\label{ccetc}
 E_\omega^0(T_1)=1+\sum_{i\le0}u_1Q^{i+1}(2,2,1,1,1,0,2,2,1)^T,\no
 \end{equation}
 where $Q$ is the mean offspring matrix of a homogeneous multitype branching process (See (\ref{cq}) below.) and $u_1$ is the mean of the immigration  (See (\ref{u1}) below.).
 But by the Ward Equation, one has that
 \begin{equation}\label{wardc}
   E_\omega^0(X_{T_1})=E_\omega^0(T_1)E_\omega^0(X_1)=\Big(1+\sum_{i\le0}u_1Q^{i+1}(2,2,1,1,1,0,2,2,1)^T\Big)(p_1+2p_2-q_1-2q_2).
 \end{equation}
 We test (with the help of Matlab) that the right-most-hand sides of  (\ref{ext1}) and (\ref{wardc}) are the same.

 At last, as application of the branching structure, we prove a law of large numbers for transient (2-2) random walk in ergodic random environment by a method known as ``the environment viewed from particles".

 We first define (2-2) random walk in random environment. Let $\Omega$ be the collection of $\omega=\{\omega_i\}_{i\in \mathbb{Z}}$ and
 $\mathcal{F}$ be the Borel $\sigma$-algebra on $\Omega.$ Define the shift operator on $\Omega$ by
 $$(\theta\omega)_i=\omega_{i+1}.$$
 Let $\mathbb{P}$ be  a probability measure on $(\Omega,\mathcal{F})$ making $(\Omega,\mathcal F,\mathbb P,\theta)$ an ergodic system.
The so-called random environment is a random element $\omega\in\Omega$ chosen according to the probability  $\mathbb{P}.$

   The (2-2) random walk $\{X_n\}_{n\ge0}$ in random environment $\omega$   is define to be  a Markov chain with initial value $X_0=x$ and transition probabilities
$$P_\omega(X_{n+1}=i+j\big |X_n=i)=\omega_i(j),\ j\in\{1,2,-1,-2\}.$$

The measure $P_\omega^x$ induced by $\{X_n\}$ on $(\mathbb{Z}^\mathbb{N},\mathcal{G}),$ with $\mathcal G$ the Borel $\sigma$-algebra,
 is called the quenched probability and  the probability $P^x$ defined on $( \mathbb Z^{\mathbb N},\mathcal G)$ by the relation
 $$P^x(B)=\int P^x_\omega(B)\mathbb P(d\omega),\ B\in \mathcal G $$
 is called the annealed probability.

 We show the following law of large numbers.
 \begin{theorem}\label{lln}
   Suppose that  $E^0(T_1)<\infty.$ Then one has that
   \begin{equation}\no
    \lim_{n\rto} \frac{X_n}{n}=V_P,
   \end{equation}
    for some $V_P.$
    Moreover  $$V_P=\frac{E_P\Big(\Pi(\omega)(2\omega_0(-2)+\omega_0(-1)+\omega_0(1)+2\omega_0(2))\Big)}{E_P(D(\omega))},$$ where $\Pi(\omega)$ and $D(\omega)$ are defined in (\ref{invd}) and (\ref{dom}) below.
\begin{remark}
  (1) The role the branching structure plays is to give  the invariant density $\Pi(\omega)$ explicitly in terms of the environment $\omega.$

  (2) In Br\'{e}mont  \cite{br09}, the author also proved a law of large numbers for (L-R) random walk in random environment, see Theorem 1.10 therein. But Br\'{e}mont \cite{br09} did not give the specific form of the velocity $V_P.$ The branching structure enables us to give the invariant density  $\Pi(\omega)$ explicitly, so that we can give the velocity $V_P$ explicitly.
 \end{remark}

 \end{theorem}
\section{Exit probability from $[a+1,b-1]$}\label{exitp}
In this section we calculate the exit probabilities of the walk from certain interval $(a,b).$

Fix $a<b.$ Let $\partial^+[a,b]=\{b,b+1\}$ and
$\partial^-[a,b]=\{a,a-1\}$ be the positive and negative boundaries
of $[a,b]$ correspondingly. For $k\in [a-1,b+1],$
$\zeta\in\partial^+[a,b]\cup\partial^-[a,b],$ define
\begin{equation}\label{epk}\no
\mcp_k(a,b,\zeta)=P_\omega^k(\text{the walk exits the interval
}[a+1,b-1]\text{ at }\zeta).\end{equation} For simplicity, we write
$\mcp_k(a,b,\zeta)$ as $\mcp_k(\zeta)$ temporarily.
Define $$M_k=\left(
    \begin{array}{ccc}
      -\frac{\omega_k(-1)+\omega_k(-2)}{\omega_k(-2)} & \frac{\omega_k(1)+\omega_k(2)}{\omega_k(-2)} & \frac{\omega_k(2)}{\omega_k(-2)} \\
     1  &  0& 0\\
      0& 1 & 0 \\
    \end{array}
  \right)$$ and set $\Pi_{r}^s=M_r\cdots M_s.$
\begin{proposition}\label{prop}
  Suppose that $\omega_i(-2)>0$ for all $i\in \mathbb{Z}.$  One has that \begin{equation}\label{pb1}\left\{\begin{array}{l}
\mcp_{b-1}(b)=\frac{\mathbf  e_1\Pi_{a+1}^{b-1}[ \mathbf e_2- \mathbf e_3]^T\Big(1+\displaystyle\sum_{l=a+1}^{b-1}\mathbf  e_1\Pi_{l}^{b-1} \mathbf e_1^T\Big)
- \mathbf e_1\Pi_{a+1}^{b-1} \mathbf e_1^T\displaystyle\sum_{l=a+1}^{b-1} \mathbf e_1
\Pi_{l}^{b-1}[\mathbf  e_2-\mathbf e_3]^T}{\mathbf e_1\Pi_{a+1}^{b-1} \mathbf e_1^T\displaystyle\sum_{l=a+1}^{b-1}\mathbf e_1
\Pi_{l}^{b-1}[\mathbf  e_1-\mathbf e_2]^T-\mathbf e_1\Pi_{a+1}^{b-1}[\mathbf e_1-\mathbf e_2]^T\Big(1+\displaystyle\sum_{l=a+1}^{b-1}\mathbf e_1\Pi_{l}^{b-1}\mathbf e_1^T\Big)},\\
\mcp_{b-2}(b)=\frac{\mathbf e_1\Pi_{a+1}^{b-1}[\mathbf e_2-\mathbf e_3]^T\displaystyle\sum_{l=a+1}^{b-1}\mathbf e_1\Pi_{l}^{b-1}[\mathbf e_1-\mathbf e_2]^T
-\mathbf e_1\Pi_{a+1}^{b-1}[\mathbf e_1-\mathbf e_2]^T\displaystyle\sum_{l=a+1}^{b-1}\mathbf e_1
\Pi_{l}^{b-1}[\mathbf e_2-\mathbf e_3]^T}{\mathbf e_1\Pi_{a+1}^{b-1}\mathbf e_1^T\displaystyle\sum_{l=a+1}^{b-1}\mathbf e_1
\Pi_{l}^{b-1}[\mathbf e_1-\mathbf e_2]^T-\mathbf e_1\Pi_{a+1}^{b-1}[\mathbf e_1-\mathbf e_2]^T\Big(1+\displaystyle\sum_{l=a+1}^{b-1}\mathbf e_1\Pi_{l}^{b-1}\mathbf e_1^T\Big)},\end{array}\right.
\end{equation}
and
\begin{equation}\label{pb2}\left\{\begin{array}{l}
\mcp_{b-1}(b+1)=\frac{\mathbf e_1\Pi_{a+1}^{b-1}\mathbf e_3^T\Big(1+\displaystyle\sum_{l=a+1}^{b-1}\mathbf e_1\Pi_{l}^{b-1}\mathbf e_1^T\Big)
-\mathbf e_1\Pi_{a+1}^{b-1}\mathbf e_1^T\displaystyle\sum_{l=a+1}^{b-1}\mathbf e_1
\Pi_{l}^{b-1}\mathbf e_3^T}{\mathbf e_1\Pi_{a+1}^{b-1}\mathbf e_1^T\displaystyle\sum_{l=a+1}^{b-1}\mathbf e_1
\Pi_{l}^{b-1}[\mathbf e_1-\mathbf e_2]^T-\mathbf e_1\Pi_{a+1}^{b-1}[\mathbf e_1-\mathbf e_2]^T\Big(1+\displaystyle\sum_{l=a+1}^{b-1}\mathbf e_1\Pi_{l}^{b-1}\mathbf e_1^T\Big)},\\
\mcp_{b-2}(b+1)=\frac{\mathbf e_1\Pi_{a+1}^{b-1}\mathbf e_3^T\displaystyle\sum_{l=a+1}^{b-1}\mathbf e_1\Pi_{l}^{b-1}[\mathbf e_1-\mathbf e_2]^T
-\mathbf e_1\Pi_{a+1}^{b-1}[\mathbf e_1-\mathbf e_2]^T\displaystyle\sum_{l=a+1}^{b-1}\mathbf e_1
\Pi_{l}^{b-1}\mathbf e_3^T}{\mathbf e_1\Pi_{a+1}^{b-1}\mathbf e_1^T\displaystyle\sum_{l=a+1}^{b-1}\mathbf e_1
\Pi_{l}^{b-1}[\mathbf e_1-\mathbf e_2]^T-\mathbf e_1\Pi_{a+1}^{b-1}[\mathbf e_1-\mathbf e_2]^T\Big(1+\displaystyle\sum_{l=a+1}^{b-1}\mathbf e_1\Pi_{l}^{b-1}\mathbf e_1^T\Big)}.\end{array}\right.
\end{equation}
\end{proposition}
\begin{remark}\label{rf}\textbf{1)} If $\omega_{i}(-2)=0$ for all $i,$ then the random walk degenerates  to (1-2) random walk. The exit probabilities could be calculated analogously as in Br\'{e}mont \cite{br02}.
\textbf{2)} If $\omega_{i}(2)=0$ for all $i\in \mathbb Z,$ then the random walk degenerates  to (2-1) random walk. For the walk transient to the right, the exit probability $\mathcal P_\omega^k(-\infty,a,a)=1$ for all $k\le a-1.$
\end{remark}

\noindent\textbf{Proof of Proposition \ref{prop}:}

 One follows from
the Markov property that
\begin{equation*}
  \mcp_k(\zeta)=\omega_k(2)\mcp_{k+2}(\zeta)+\omega_{k}(1)\mcp_{k+1}(\zeta)+\omega_{k}(-1)\mcp_{k-1}(\zeta)+\omega_{k}(-2)\mcp_{k-2}(\zeta),
\end{equation*}
which leads to
\begin{equation}\label{hom}\begin{split}
  (\mcp_{k-1}-\mcp_{k-2})(\zeta)=&\frac{\omega_k(-1)+\omega_k(-2)}{\omega_k(-2)}(\mcp_{k}-\mcp_{k-1})(\zeta)\\
  &+\frac{\omega_k(1)+\omega_k(2)}{\omega_k(-2)}(\mcp_{k+1}-\mcp_{k})(\zeta)+\frac{\omega_k(2)}{\omega_k(-2)}(\mcp_{k+2}-\mcp_{k+1})(\zeta).
\end{split}\end{equation}
Writing the equation (\ref{hom}) in the matrix form, one has that
\begin{equation}\label{mform}
  \left(
    \begin{array}{ccc}
     - \frac{\omega_k(-1)+\omega_k(-2)}{\omega_k(-2)} & \frac{\omega_k(1)+\omega_k(2)}{\omega_k(-2)} & \frac{\omega_k(2)}{\omega_k(-2)} \\
     1  &  0& 0\\
      0& 1 & 0 \\
    \end{array}
  \right)\left(
           \begin{array}{c}
             (\mcp_k-\mcp_{k-1})(\zeta) \\
             (\mcp_{k+1}-\mcp_{k})(\zeta) \\
             (\mcp_{k+2}-\mcp_{k+1})(\zeta) \\
           \end{array}
         \right)=\left(\begin{array}{c}
             (\mcp_{k-1}-\mcp_{k-2})(\zeta) \\
             (\mcp_{k}-\mcp_{k-1})(\zeta) \\
             (\mcp_{k+1}-\mcp_{k})(\zeta) \\
           \end{array}
         \right).
\end{equation}
Define
\begin{equation}\label{mk}\no
  V_k(\zeta)=\left(\begin{array}{c}
             (\mcp_{k-1}-\mcp_{k-2})(\zeta) \\
             (\mcp_{k}-\mcp_{k-1})(\zeta) \\
             (\mcp_{k+1}-\mcp_{k})(\zeta) \\
           \end{array}
         \right).
\end{equation}
Then equation (\ref{mform}) becomes
\begin{equation}\label{ditui}
  V_k(\zeta)=M_kV_{k+1}(\zeta).
\end{equation}
We note that the equation  $ V_k(\zeta)=M_kV_{k+1}$ makes sense for $k\in(a,b].$ Since $\mcp_{a-2}(a,b,\zeta)$ has no sense, so is $V_a(\zeta).$
Since $\mcp_b(a,b,b)=1,\mcp_b(a,b,b+1)=0,$ and $\mcp_{b+1}(a,b,b)=0,\mcp_{b+1}(a,b,b+1)=1,$ then
\begin{equation}\label{boundb}\no
  V_{b}(b)=\left(\begin{array}{c}
             (\mcp_{b-1}-\mcp_{b-2})(b) \\
             (\mcp_{b}-\mcp_{b-1})(b) \\
             (\mcp_{b+1}-\mcp_{b})(b) \\
           \end{array}
         \right)=\left(\begin{array}{c}
             (\mcp_{b-1}-\mcp_{b-2})(b) \\
             1-\mcp_{b-1}(b) \\
             -1 \\
           \end{array}
         \right)
\end{equation} and
\begin{equation}\label{boundb1}\no
  V_b({b+1})=\left(\begin{array}{c}
             (\mcp_{b-1}-\mcp_{b-2})(b+1) \\
             (\mcp_{b}-\mcp_{b-1})(b+1) \\
             (\mcp_{b+1}-\mcp_{b})(b+1) \\
           \end{array}
         \right)=\left(\begin{array}{c}
             (\mcp_{b-1}-\mcp_{b-2})(b+1) \\
             -\mcp_{b-1}(b+1) \\
             1 \\
           \end{array}
         \right).
\end{equation}
Substituting to (\ref{ditui}) one has that
$$\mathbf e_1V_k(b)=\z(\mcp_{k-1}-\mcp_{k-2}\y)(b)=\mathbf e_1M_k\cdots M_{b-1}\z(\mcp_{b-1}(b)[\mathbf e_1-\mathbf e_2]-\mcp_{b-2}(b)\mathbf e_1+[\mathbf e_2-\mathbf e_3]\y)^T.$$
Summing from $k$ to $b-1,$ it follows that
\begin{equation}\label{pkb}
\z(\mcp_{b-2}-\mcp_{k-2}\y)(b)=\sum_{l=k}^{b-1}\mathbf e_1M_l\cdots M_{b-1}\Big(\mcp_{b-1}(b)[\mathbf e_1-\mathbf e_2]
-\mcp_{b-2}(b)\mathbf e_1+[\mathbf e_2-\mathbf e_3]\Big)^T.\end{equation}
Note that $\mcp_{a}(a,b,b)=0.$ Then
$$\mcp_{b-2}(b)=\sum_{l=a+2}^{b-1}\mathbf e_1M_l\cdots M_{b-1}\Big(\mcp_{b-1}(b)[\mathbf e_1-\mathbf e_2]-\mcp_{b-2}(b)\mathbf e_1+[\mathbf e_2-\mathbf e_3]\Big)^T.$$
On the other hand, since
\begin{equation}\no
  0=(\mcp_{a}-\mcp_{a-1})(b)=\mathbf e_1M_{a+1}\cdots M_{b-1}\Big(\mcp_{b-1}(b)[\mathbf e_1-\mathbf e_2]-\mcp_{b-2}(b)\mathbf e_1+[\mathbf e_2-\mathbf e_3]\Big)^T
\end{equation}
then
\begin{equation}\label{eqs}\left\{\begin{array}{l}
  \mathbf e_1M_{a+1}\cdots M_{b-1}\Big(\mcp_{b-1}(b)[\mathbf e_1-\mathbf e_2]-\mcp_{b-2}(b)\mathbf e_1+[\mathbf e_2-\mathbf e_3]\Big)^T=0,\\
  \mcp_{b-2}(b)=\displaystyle\sum_{l=a+1}^{b-1}\mathbf e_1M_l\cdots M_{b-1}\Big(\mcp_{b-1}(b)[\mathbf e_1-\mathbf e_2]-\mcp_{b-2}(b)\mathbf e_1+[\mathbf e_2-\mathbf e_3]\Big)^T.\end{array}\right.
\end{equation}
Solving (\ref{eqs}), one gets (\ref{pb1}).

Also, for $a+1\le k\le b-3,$ one has from (\ref{pkb}) that
\begin{equation}\label{pkbs}\no\begin{split}
\mcp_{k}(b)=&\mcp_{b-2}(b)\Big(1+\sum_{l=k+2}^{b-1}\mathbf e_1\Pi_{l}^{b-1}\mathbf e_1^T\Big)\\
&-\mcp_{b-1}(b)
\sum_{l=k+2}^{b-1}\mathbf e_1\Pi_{l}^{b-1}[\mathbf e_1-\mathbf e_2]^T-\sum_{l=k+2}^{b-1}\mathbf e_1\Pi_{l}^{b-1}[\mathbf e_2-\mathbf e_3]^T.
\end{split}
\end{equation}

Next we calculate $\mcp_{b-1}(b+1)$ and $\mcp_{b-2}(b+1).$ Similarly as (\ref{eqs}), one has that
\begin{equation}\label{eqsb}\no\left\{\begin{array}{l}
  \mathbf e_1M_{a+1}\cdots M_{b-1}\Big(\mcp_{b-1}(b+1)[\mathbf e_1-\mathbf e_2]-\mcp_{b-2}(b+1)\mathbf e_1+ \mathbf e_3\Big)^T=0,\\
  \mcp_{b-2}(b+1)=\displaystyle\sum_{l=a+1}^{b-1}\mathbf e_1M_l\cdots M_{b-1}\Big(\mcp_{b-1}(b+1)[\mathbf e_1-\mathbf e_2]-\mcp_{b-2}(b+1)\mathbf e_1+\mathbf e_3\Big)^T\end{array}\right.
\end{equation}
which leads to (\ref{pb2}).

Also for $a+1\le k\le b-3,$ one has that

\begin{equation}\label{pkbs1}\no\begin{split}
\mcp_{k}(b+1)=&\mcp_{b-2}(b+1)\Big(1+\sum_{l=k+2}^{b-1}\mathbf e_1\Pi_{l}^{b-1}\mathbf e_1^T\Big)\\
&-\mcp_{b-1}(b+1)
\sum_{l=k+2}^{b-1}\mathbf e_1\Pi_{l}^{b-1}[\mathbf e_1-\mathbf e_2]^T-\sum_{l=k+2}^{b-1}\mathbf e_1\Pi_{l}^{b-1}\mathbf e_3^T.
\end{split}
\end{equation}
\qed

\section{Path decomposition--Proofs of Theorem \ref{thbranc}}\label{bran}

Throughout this section, we always assume that $\limsup_{n\rto} X_n=\infty.$ That is, the walk $\{X_n\}$ is transient to $\infty$ or recurrent. The notation $\{i\rightarrow j\}$ will  be always used to denote a jump (a step) by the walk from $i$ to $j.$

Define $$T_1=\inf\{n\ge0:X_n>0\},$$ the hitting time of
$[1,\infty).$ The purpose of this section is to count exactly all
steps by the walk before $T_1.$
\subsection{The excursions, corresponding probabilities and immigration distributions}
From Proposition \ref{prop} and Remark \ref{rf}, we see that the exit probabilities of the walk from certain interval $(a,b)$ could be expressed in terms of $\omega.$ So we always assume that all these exit probabilities are already known in the remainder of the paper.

For $k\le i<j,$ denote
\begin{equation}\no
  f_k(i,j)=P_\omega^k(\text{the walk hits } (i,\infty) \text{ from below at }j).
\end{equation}
We remark that for (2-2) random walk, in the definition of $f_k(i,j),$  the term $j$ only takes values in $\{i+1,i+2\}.$ With the notations of Section  \ref{exitp},
$$f_k(i,i+1)=\mcp_k(-\infty,i+1,i+1)\text{ and } f_k(i,i+2)=\mcp_k(-\infty,i+1,i+2).$$
Next we analyze the path of the walk.

Firstly, we consider a special excursion, which will be called
type-$\mathcal{A}$ excursion latter, of the walk.
\begin{definition}
We call excursions of the form $\{X_k=i,X_{k+1}=i-1,X_{k+2}\le
i-1,...,X_{k+l}\le i-1,X_{k+l+1}\ge i\}$ type-$\mathcal{A}$
excursions at $i.$ Corresponding to the three kinds of possible last
step of type-$\mathcal{A}$ excursions $i,$ say, $\{i-1\rightarrow i\},$  $\{i-2\rightarrow i\}$ and
$\{i-1\rightarrow i+1\},$ we classify type-$\mathcal{A}$ excursions at $i$ into three
sub-types $\mathcal A_{i,1},$ $\mathcal A_{i,2}$ and $\mathcal A_{i,3}.$
\end{definition}

An excursion will be also called a particle some times in the remainder of the
paper.

 Note that a type-$\mathcal{A}$ particle at $i$ begins when the walk jumps down from $i$ to $i-1.$ After
that the walk runs in $(-\infty,i-1].$ At last the walk hits
$[i,\infty)$ at some $j$ and the excursion goes to end.

 Next we define some indexes $\alpha_{i,1},$ $\alpha_{i,3}$  and
$\alpha_{i,2}$ correspondingly to $\mathcal A_{i,1}$ $\mathcal A_{i,3},$ and
$\mathcal A_{i,2}.$ Let
\begin{equation}\no
  \begin{split}
&\alpha_{i,1}:=\omega_{i}(-1)\sum_{n,m\ge0}\frac{(n+m)!}{n!m!}[\omega_{i-1}(-1)f_{i-2}(i-2,i-1)]^n[\omega_{i-1}(-2)f_{i-3}(i-2,i-1)]^m\omega_{i-1}(1),\\
&\alpha_{i,3}:=\omega_{i}(-1)\sum_{n,m\ge0}\frac{(n+m)!}{n!m!}[\omega_{i-1}(-1)f_{i-2}(i-2,i-1)]^n[\omega_{i-1}(-2)f_{i-3}(i-2,i-1)]^m\omega_{i-1}(2),\\
&\alpha_{i,2}:=\omega_{i}(-1)-\alpha_{i,1}-\alpha_{i,3}.
  \end{split}
\end{equation}
Note that  $\alpha_{i,1}$ differs from $\alpha_{i,3}$ only in the
last term of the product. Therefore we explain only the meaning of
$\alpha_{i,1}.$ The first term $\omega_{i}(-1)$ is transition
probability of the first step of the excursion $\mathcal A_{i,1}$ from $i$ to
$i-1.$ The last term $\omega_{i-1}(1)$ is transition probability
 of the last step of the excursion $\mathcal A_{i,1}$ from $i-1$ to $i.$
 The summation in the center indicates all events occurring between the
 first step and the last step. In details, before the last step happens,
 $n$ steps of the form $\{i-1\rightarrow i-2\}$ and $m$ steps of the form $\{i-1\rightarrow i-3\}$
 occur and the total number of possible combinations of these $m+n$ steps is
 $\left(
                                                                                 \begin{array}{c}
                                                                                   n+m \\
                                                                                   m\\
                                                                                 \end{array}
                                                                               \right)
=\frac{(n+m)!}{n!m!}.$
The term $\omega_{i-1}(-1)f_{i-2}(i-2,i-1)$ means that, with probability $\omega_{i-1}(-1),$ a step $\{i-1\rightarrow i-2\}$ occurs,
and  starting from $i-2,$ with probability $f_{i-2}(i-2,i-1),$ it hits $[i-1,\infty)$ at $i-1.$
The term $\omega_{i-1}(-2)f_{i-3}(i-2,i-1)$ could be explained analogously.

In fact, $\alpha_{i,2}$ could be defined similarly as $\alpha_{i,1}$
and $\alpha_{i,3}.$ But the definition is tedious. We note that after the
walk jumps down with probability
$\omega_{i}(-1)$ from $i$ to $i-1,$  then starting from $i-1,$ it hits $[i,\infty)$
with probability $1.$ Then the summation of $\alpha_{i,1},$
$\alpha_{i,2},$ and $\alpha_{i,3}$ should be $\omega_i({-1}).$ So we
define
$$\alpha_{i,2}:=\omega_{i}(-1)-\alpha_{i,1}-\alpha_{i,3}.$$

Some easy calculation shows that
\begin{equation}\no
  \begin{split}
&\alpha_{i,1}=\frac{\omega_{i}(-1)\omega_{i-1}(1)}{1-\omega_{i-1}(-1)f_{i-2}(i-2,i-1)-\omega_{i-1}(-2)f_{i-3}(i-2,i-1)},\\
&\alpha_{i,3}=\frac{\omega_{i}(-1)\omega_{i-1}(2)}{1-\omega_{i-1}(-1)f_{i-2}(i-2,i-1)-\omega_{i-1}(-2)f_{i-3}(i-2,i-1)},\\
&\alpha_{i,2}=\frac{\omega_{i}(-1)[1-\omega_{i-1}(1)-\omega_{i-1}(2)-\omega_{i-1}(-1)f_{i-2}(i-2,i-1)-\omega_{i-1}(-2)f_{i-3}(i-2,i-1)]}{1-\omega_{i-1}(-1)f_{i-2}(i-2,i-1)-\omega_{i-1}(-2)f_{i-3}(i-2,i-1)}.
  \end{split}
\end{equation}
Secondly, we describe another excursion, that is, type-$\mathcal{B}$ excursion.
\begin{definition}
We call excursions of the form $\{X_k=i,X_{k+1}=i-2,X_{k+2}\le
i-1,...,X_{k+l}\le i-1,X_{k+l+1}\ge i\}$ type-$\mathcal{B}$
excursions at $i.$ Corresponding to the three kinds of possible last
step of type-$\mathcal{B}$ excursions at $i,$ say, $\{i-1\rightarrow i\},$ $\{i-2\rightarrow i\}$ and
$\{i-1\rightarrow i+1\},$ we classify type-$\mathcal{B}$ excursions at $i$ into three
sub-types $\mathcal B_{i,1},$ $\mathcal B_{i,2}$ and $\mathcal B_{i,3}.$
\end{definition}
Note that a type-$\mathcal{B}$ excursion begins when the walk jumps down from
$i$ to $i-2.$ After that the walk runs in $(-\infty,i-1].$ At last
the walk hits $[i,\infty)$ at some $j$ and the excursion goes to
end.

 Next we define some indexes $\beta_{i,1}$ $\beta_{i,2}$ and $\beta_{i,3}$ corresponding to $\mathcal B_{i,1}$ $\mathcal B_{i,3},$ and $\mathcal B_{i,2}.$
Let
\begin{equation}\label{indb}
  \begin{split}
&\beta_{i,1}:=\omega_{i}(-2)f_{i-2}(i-2,i-1)\\
&\hspace{1.5cm}\times\sum_{n,m\ge0}\frac{(n+m)!}{n!m!}[\omega_{i-1}(-1)f_{i-2}(i-2,i-1)]^n[\omega_{i-1}(-2)f_{i-3}(i-2,i-1)]^m\omega_{i-1}(1),\\
&\beta_{i,3}:=\omega_{i}(-2)f_{i-2}(i-2,i-1)\\
&\hspace{1.5cm}\times\sum_{n,m\ge0}\frac{(n+m)!}{n!m!}[\omega_{i-1}(-1)f_{i-2}(i-2,i-1)]^n[\omega_{i-1}(-2)f_{i-3}(i-2,i-1)]^m\omega_{i-1}(2),\\
&\beta_{i,2}:=\omega_{i}(-2)-\beta_{i,1}-\beta_{i,3}.
  \end{split}
\end{equation}

$\beta_{i,1}$ differs from $\alpha_{i,1}$ only in the  term $\omega_{i}(2)f_{i-2}(i-2,i-1).$ We explain only this term.
Starting from $i,$ with probability $\omega_i(-2),$ the walk jumps down from $i$ to $i-2$  and the excursion begins. Recall that $\beta_{i,1}$ is the index corresponding to a
$\mathcal B_{i,1}$ excursion. Since a $\mathcal B_{i,1}$ excursion ends with a jump $\{i-1\rightarrow i\},$ after visiting $i-2,$ it will reach $i-1$ from below before it comes to end. The sum of products of transition probabilities of all possible paths from $i-2$ to hit $i-1$ from below is $f_{i-2}(i-2,i-1).$

One follows from (\ref{indb}) that
\begin{eqnarray*}
&&\beta_{i,1}=\frac{\omega_{i}(-2)f_{i-2}(i-2,i-1)\omega_{i-1}(1)}{1-\omega_{i-1}(-1)f_{i-2}(i-2,i-1)-\omega_{i-1}(-2)f_{i-3}(i-2,i-1)},\\
&&\beta_{i,3}=\frac{\omega_{i}(-2)f_{i-2}(i-2,i-1)\omega_{i-1}(2)}{1-\omega_{i-1}(-1)f_{i-2}(i-2,i-1)-\omega_{i-1}(-2)f_{i-3}(i-2,i-1)},\\
&&\beta_{i,2}=\frac{\omega_{i}(-2)[1-(\omega_{i-1}(-1)+\omega_{i-1}(1)+\omega_{i-1}(2))f_{i-2}(i-2,i-1)-\omega_{i-1}(-2)f_{i-3}(i-2,i-1)]}{1-\omega_{i-1}(-1)f_{i-2}(i-2,i-1)-\omega_{i-1}(-2)f_{i-3}(i-2,i-1)}.
  \end{eqnarray*}

At last, we define the third kind of excursion, type-$\mathcal{C}$ excursion.
\begin{definition}
We call excursions of the form $\{X_k=i+1,X_{k+1}=i-1,X_{k+2}\le
i-1,...,X_{k+l}\le i-1,X_{k+l+1}\ge i\}$ type-$\mathcal{C}$
excursions at $i.$ Corresponding to the three kinds of possible last
step of type-$\mathcal{C}$ excursions $i,$ say, $\{i-1\rightarrow i\},$ $\{i-2\rightarrow i\}$ and
$\{i-1\rightarrow i+1\},$  we classify type-$\mathcal{C}$ excursions at $i$ into three
sub-types $\mathcal C_{i,1},$ $\mathcal C_{i,2}$ and $\mathcal C_{i,3}.$
\end{definition}
Note that a type-$\mathcal{C}$ excursion at $i$ begins when the walk jumps down from
$i+1$ to $i-1.$ After that the walk runs in $(-\infty,i-1].$ At last
the walk hits $[i,\infty)$ at some $j$ and the excursion goes to
end.

 Next we define some indexes  $\gamma_{i,1}$ $\gamma_{i,2}$ and $\gamma_{i,3}$ corresponding to $\mathcal C_{i,1}$ $\mathcal C_{i,3},$ and $\mathcal C_{i,2}.$
Let
\begin{equation}\label{indc}
  \begin{split}
&\gamma_{i,1}:=\omega_{i+1}(-2)\sum_{n,m\ge0}\frac{(n+m)!}{n!m!}[\omega_{i-1}(-1)f_{i-2}(i-2,i-1)]^n[\omega_{i-1}(-2)f_{i-3}(i-2,i-1)]^m\omega_{i-1}(1),\\
&\gamma_{i,3}:=\omega_{i+1}(-2)\sum_{n,m\ge0}\frac{(n+m)!}{n!m!}[\omega_{i-1}(-1)f_{i-2}(i-2,i-1)]^n[\omega_{i-1}(-2)f_{i-3}(i-2,i-1)]^m\omega_{i-1}(2),\\
&\gamma_{i,2}:=\omega_{i+1}(-2)-\gamma_{i,1}-\gamma_{i,3}.
  \end{split}
\end{equation}
A $\mathcal C_{i,1}$ excursion differs from an $\mathcal A_{i,1}$ excursion only in the first step. The first term in the
 product of $\gamma_{i,1}$ is $\omega_{i+1}(-2).$ It means that with probability $\omega_{i+1}(-2),$
  the walk jumps down from $i+1$ to $i-1,$ and the excursion begins.
   $\gamma_{i,2}$ and $\gamma_{i,3}$ could be understood analogously as $\alpha_{i,2}$
   and $\alpha_{i,3}.$ We will not repeat them here.

   One follows from (\ref{indc}) that
   \begin{equation}\no
  \begin{split}
&\gamma_{i,1}=\frac{\omega_{i+1}(-2)\omega_{i-1}(1)}{1-\omega_{i-1}(-1)f_{i-2}(i-2,i-1)-\omega_{i-1}(-2)f_{i-3}(i-2,i-1)},\\
&\gamma_{i,3}=\frac{\omega_{i+1}(-2)\omega_{i-1}(2)}{1-\omega_{i-1}(-1)f_{i-2}(i-2,i-1)-\omega_{i-1}(-2)f_{i-3}(i-2,i-1)},\\
&\gamma_{i,2}=\frac{\omega_{i+1}(-2)[1-\omega_{i-1}(1)-\omega_{i-1}(2)-\omega_{i-1}(-1)f_{i-2}(i-2,i-1)-\omega_{i-1}(-2)f_{i-3}(i-2,i-1)]}{1-\omega_{i-1}(-1)f_{i-2}(i-2,i-1)-\omega_{i-1}(-2)f_{i-3}(i-2,i-1)}.
  \end{split}
\end{equation}

Define $$A_{i,j}=\#\{\mathcal A_{i,j} \text{ excursions before } T_1\},$$
$$B_{i,j}=\#\{\mathcal B_{i,j} \text{ excursions before } T_1\},$$
$$C_{i,j}=\#\{\mathcal C_{i,j} \text{ excursions before } T_1\},$$ for $i\le 0$ and $j=1,2,3.$

We aim at counting exactly all steps by the walk before $T_1.$ For this purpose, define
\begin{equation}\no
  U_{i}=(A_{i,1},A_{i,2},A_{i,3},B_{i,1},B_{i,2},B_{i,3},C_{i,1},C_{i,2},C_{i,3})
\end{equation} being the  total number of different excursions at $i$ before time $T_1.$
Next we show by path decomposition that $\{U_i\}_{i\le0}$ is a non-homogeneous
multitype branching process.
\begin{center}
\includegraphics[height=6cm]{immigration.1}
\end{center}

Firstly, the branching process needs an ancestor (some particle immigrating in). The walk starts from
$0.$ But before $T_1,$ there is no jump down from above $1$ to $0$
by the walk. One can imagine that there is a step by the walk from
$1$ to $0$ before it starts from $0$ (One can also imagine that this
step is from $2$ to $0.$ But this makes no difference.), that is,
set $X_{-1}=1.$ Adding this imaginary step, the path $\{X_{-1}=1,
X_0=0,X_1,...,X_{T_1}\}$ forms a type-$\mathcal{A}$ excursion at $1$  such that
$$A_{1,1}+A_{1,2}+A_{1,3}=1.$$

The distributions of $A_{1,1},$
$A_{1,2},$ and $A_{1,3}$ are
\begin{equation*}
  \begin{split}
    P_\omega^0(A_{1,1}=1)
    =\frac{\alpha_{1,1}}{\alpha_{1,1}+\alpha_{1,2}+\alpha_{1,3}}
   =\frac{\omega_{0}(1)}{1-\omega_{0}(-1)f_{-1}(-1,0)-\omega_{0}(-2)f_{-2}(-1,0)},\end{split}
  \end{equation*}
 \begin{equation*}
  \begin{split}
    P_\omega^0(A_{1,3}=1)
    =\frac{\alpha_{1,3}}{\alpha_{1,1}+\alpha_{1,2}+\alpha_{1,3}}
   =\frac{\omega_{0}(2)}{1-\omega_{0}(-1)f_{-1}(-1,0)-\omega_{0}(-2)f_{-2}(-1,0)},\end{split}
  \end{equation*}
  and
  \begin{equation*}
  \begin{split}
   P_\omega^0(A_{1,2}=1&)=\frac{\alpha_{1,2}}{\alpha_{1,1}+\alpha_{1,2}+\alpha_{1,3}}\\
   &=\frac{[1-\omega_{0}(1)-\omega_{0}(2)-\omega_{0}(-1)f_{-1}(-1,0)-\omega_{0}(-2)
  f_{-2}(-1,0)]}{1-\omega_{0}(-1)f_{-1}(-1,0)-\omega_{0}(-2)f_{-2}(-1,0)}.
   \end{split}
  \end{equation*}
  The meaning of $P_\omega^0(A_{1,1}=1)$ is obvious. The total sum of
  the product of transition probabilities of all excursions of the form
  $\{X_{-1}=0,X_{0}=1,X_1,...,X_{T_1}\}$ is
  $\alpha_{1,1}+\alpha_{1,2}+\alpha_{1,3}=\omega_1(-1)$ and the total sum of
  the product of transition probabilities of all possible paths of an $\mathcal A_{1,1}$ excursion
  is
  $\alpha_{1,1}$. Therefore $P_\omega^0(A_{1,1}=1)
    =\frac{\alpha_{1,1}}{\alpha_{1,1}+\alpha_{1,2}+\alpha_{1,3}}.$
 The values of $P_\omega^0(A_{1,2}=1)$ and $ P_\omega^0(A_{1,3}=1)$
 could be explained analogously.

We can treat the above discussed imaginary excursion as the particle
immigrates in the branching system, say, the ancestor (immigration) of the
branching process. The immigration laws are
\begin{equation}\label{inia}
  \begin{split}
    P_\omega^0(U_1=(1,0,...,0))
    =\frac{\omega_{0}(1)}{1-\omega_{0}(-1)f_{-1}(-1,0)-\omega_{0}(-2)f_{-2}(-1,0)},\end{split}
  \end{equation}
 \begin{equation}\label{inib}
  \begin{split}
    P_\omega^0(U_1=(0,1,0,...,0))
       =\frac{\omega_{0}(2)}{1-\omega_{0}(-1)f_{-1}(-1,0)-\omega_{0}(-2)f_{-2}(-1,0)},\end{split}
  \end{equation}
  and
  \begin{equation}\label{inic}
  \begin{split}
   P_\omega^0(U_1=(0,0,1,0,...,0))
   =\frac{[1-\omega_{0}(1)-\omega_{0}(2)-\omega_{0}(-1)f_{-1}(-1,0)-\omega_{0}(-2)
  f_{-2}(-1,0)]}{1-\omega_{0}(-1)f_{-1}(-1,0)-\omega_{0}(-2)f_{-2}(-1,0)}.
   \end{split}
  \end{equation}
\subsection{Branching mechanisms}
  After revealing the immigration, we discuss the branching mechanism. Although there are 9 types of particles in the system, many of them share the same offspring distributions.

{\bf \noindent(a) Offspring distributions of $\mathcal A_{i+1,1},$  $\mathcal A_{i+1,3},$  $\mathcal C_{i+1,1}$  and $\mathcal C_{i+1,3}$ particles}

For $i\le 0,$ conditioned on $U_{i+1}=(1,0,...,0),$ that is,
$\{A_{i+1,1}=1,A_{i+1,2}=A_{i+1,3}=B_{i+1,1}=B_{i+1,2}=B_{i+1,3}=C_{i+1,1}=C_{i+1,2}=C_{i+1,3}=0\},$
we study the distribution of $U_{i}.$ Note that the first step of
excursion $\mathcal A_{i+1,1}$ is $\{i+1\rightarrow i\}$ and the last step is
$\{i\rightarrow i+1\}.$ All contributions of this excursion to $U_i$ occurs
between these two steps. Therefore this particle could give
birthes only to excursions $\mathcal A_{i,1},$ $\mathcal A_{i,2},$ $\mathcal B_{i,1},$ and $\mathcal B_{i,2}.$ We
calculate the probability of the event
$$\{A_{i,1}=a, A_{i,2}=b, B_{i,1}=c, B_{i,2}=d\},$$ that is,
$$\{U_{i}=(a,b,0,c,d,0,0,0,0)\},$$
conditioned on $\{A_{i+1,1}=1\}$(Or $U_{i+1}=(1,0,...,0)$). Indeed, one follows from the strong Markov property that these $a+b+c+d$ excursions at $i$ are independent,
and the total number of all possible combinations of those excursions is $\frac{(a+b+c+d)!}{a!b!c!d!}.$ Therefore
\begin{equation}\label{oaa}
\begin{split}
  P_\omega^0(U_{i}=(a,b,&0,c,d,0,0,0,0)\big|U_{i+1}=(1,0,...,0))\\
  &=\frac{(a+b+c+d)!}{a!b!c!d!}\alpha_{i,1}^a\alpha_{i,2}^b\beta_{i,1}^c\beta_{i,2}^d(1-\alpha_{i,1}-\alpha_{i,2}-\beta_{i,1}-\beta_{i,2}).
\end{split}
\end{equation}
\begin{center}\includegraphics[height=6cm]{offspringa.1}\end{center}

One notes that
$\mathcal A_{i+1,3},$ $\mathcal C_{i+1,1},$  $\mathcal C_{i+1,3}$ and $\mathcal A_{i+1,1}$ excursions share a common property, that is, the first step is from above $i$ to $i$ and the last step is from $i$ to above $i.$ Therefore they share the same offspring distribution. So analogously, one has that
\begin{equation}\label{oab}
\begin{split}
  P_\omega^0(U_{i}=(a,b,&0,c,d,0,0,0,0)\big|U_{i+1}=(0,0,1,0,...,0))\\
  &=\frac{(a+b+c+d)!}{a!b!c!d!}\alpha_{i,1}^a\alpha_{i,2}^b\beta_{i,1}^c\beta_{i,2}^d(1-\alpha_{i,1}-\alpha_{i,2}-\beta_{i,1}-\beta_{i,2}),
\end{split}
\end{equation}
\begin{equation}\label{oac}
\begin{split}
  P_\omega^0(U_{i}=(a,b,&0,c,d,0,0,0,0)\big|U_{i+1}=(0,..,0,1,0,0))\\
  &=\frac{(a+b+c+d)!}{a!b!c!d!}\alpha_{i,1}^a\alpha_{i,2}^b\beta_{i,1}^c\beta_{i,2}^d(1-\alpha_{i,1}-\alpha_{i,2}-\beta_{i,1}-\beta_{i,2}),
\end{split}
\end{equation}
\begin{equation}\label{oad}
\begin{split}
  P_\omega^0(U_{i}=(a,b,&0,c,d,0,0,0,0)\big|U_{i+1}=(0,..,0,1))\\
  &=\frac{(a+b+c+d)!}{a!b!c!d!}\alpha_{i,1}^a\alpha_{i,2}^b\beta_{i,1}^c\beta_{i,2}^d(1-\alpha_{i,1}-\alpha_{i,2}-\beta_{i,1}-\beta_{i,2}).
\end{split}
\end{equation}

{\bf \noindent(b) Offspring distributions of $\mathcal A_{i+1,2},$ and $\mathcal C_{i+1,2}$ particles}

\begin{center}
  \includegraphics[height=6cm]{offspringb.1}
\end{center}

Conditioned on $\{U_{i+1}=(0,1,0,...,0)\},$ that is, $\{A_{i+1,2}=1,A_{i+1,1}=A_{i+1,3}=B_{i+1,1}=B_{i+1,2}=B_{i+1,3}=C_{i+1,1}=C_{i+1,2}=C_{i+1,3}=0\},$ we discussion the distribution of $U_i.$
Recall that the first step of an $\mathcal A_{i+1,2}$ excursion is $\{i+1\rightarrow i\}$ and the last step is $\{i-1\rightarrow i+1\}.$ Things get delicate because the last step  $\{i-1\rightarrow i+1\}.$
Before the last step occurs, the walk must jump down from $i,$ possibly to $i-1$ or $i-2.$ If it jumps down from $i$ to $i-1,$ it gives birth to an $\mathcal A_{i,3}$ particle; if it jumps down from $i$ to $i-2,$ it gives birth to a $\mathcal B_{i,3}$ particle. That is
$$P_\omega^0(A_{i,3}+B_{i,3}=1\big|U_{i+1}=\mathbf e_2)=1.$$

The sum of the product of the transition probabilities of all possible paths of an excursion $\mathcal A_{i,3}$ is $\alpha_{i,3},$ and that of a $\mathcal B_{i,3}$ excursion is $\beta_{i,3}.$
In this point of view, one concludes that
\begin{equation}\label{unab}
  P_{\omega}^0(A_{i,3}=1\big|U_{i+1}=\mathbf e_2)=1-P_{\omega}^0(B_{i,3}=1\big|U_{i+1}=\mathbf e_2)=\frac{\alpha_{i,3}}{\alpha_{i,3}+\beta_{i,3}}.
\end{equation}

Before giving birth to the above discussed particle $\mathcal B_{i,3}$ or $\mathcal A_{i,3},$ the excursion $\mathcal A_{i+1,2}$ may give birthes to a number of
$\mathcal A_{i,1},$ $\mathcal A_{i,2},$ $\mathcal B_{i,1}$ and $\mathcal B_{i,2}$ particles. Once again, one follows by path decomposition and Markov property that, all excursions born to particle $\mathcal A_{i+1,2}$ are independent. Therefore one has that
 \begin{equation}\label{oba}
\begin{split}
  P_\omega^0(U_{i}=(a,b,&1,c,d,0,0,0,0)\big|U_{i+1}=(0,1,0,...,0))\\
  &=\frac{(a+b+c+d)!}{a!b!c!d!}\alpha_{i,1}^a\alpha_{i,2}^b\beta_{i,1}^c\beta_{i,2}^d(1-\alpha_{i,1}-\alpha_{i,2}-\beta_{i,1}-\beta_{i,2})\frac{\alpha_{i,3}}{\alpha_{i,3}+\beta_{i,3}},
\end{split}
\end{equation}
 and
  \begin{equation}\label{obb}
\begin{split}
  P_\omega^0(U_{i}=(a,b,&0,c,d,1,0,0,0)\big|U_{i+1}=(0,1,0,...,0))\\
  &=\frac{(a+b+c+d)!}{a!b!c!d!}\alpha_{i,1}^a\alpha_{i,2}^b\beta_{i,1}^c\beta_{i,2}^d(1-\alpha_{i,1}-\alpha_{i,2}-\beta_{i,1}-\beta_{i,2})\frac{\beta_{i,3}}{\alpha_{i,3}+\beta_{i,3}}.
\end{split}
\end{equation}
Also, analogously, the offspring distributions of a $\mathcal C_{i+1,2}$ particle are
\begin{equation}\label{obc}
\begin{split}
  P_\omega^0(U_{i}=(a,b,&1,c,d,0,0,0,0)\big|U_{i+1}=(0,...,0,1,0))\\
  &=\frac{(a+b+c+d)!}{a!b!c!d!}\alpha_{i,1}^a\alpha_{i,2}^b\beta_{i,1}^c\beta_{i,2}^d(1-\alpha_{i,1}-\alpha_{i,2}-\beta_{i,1}-\beta_{i,2})\frac{\alpha_{i,3}}{\alpha_{i,3}+\beta_{i,3}},
\end{split}
\end{equation}
 and
  \begin{equation}\label{obd}
\begin{split}
  P_\omega^0(U_{i}=(a,b,&0,c,d,1,0,0,0)\big|U_{i+1}=(0,...,0,1,0))\\
  &=\frac{(a+b+c+d)!}{a!b!c!d!}\alpha_{i,1}^a\alpha_{i,2}^b\beta_{i,1}^c\beta_{i,2}^d(1-\alpha_{i,1}-\alpha_{i,2}-\beta_{i,1}-\beta_{i,2})\frac{\beta_{i,3}}{\alpha_{i,3}+\beta_{i,3}}.
\end{split}
\end{equation}

{\bf\noindent (C) Offspring distributions of $\mathcal B_{i+1,1},$ and $\mathcal B_{i+1,3}$ particles }
\begin{center}
  \includegraphics[height=7cm]{offspringc.1}
\end{center}

$\mathcal B_{i+1,1}$ and $\mathcal B_{i+1,3}$ excursions differ from each other only
in the last step. But their last steps are both from $i$ to above $i.$ Therefore they have the same offspring distributions.
Conditioned on $U_{i+1}=(0,0,0,1,0,...,0),$ that is,
$\{B_{i+1,1}=1,A_{i+1,1}=A_{i+1,2}=A_{i+1,3}=B_{i+1,2}=B_{i+1,3}=C_{i+1,1}=C_{i+1,2}=C_{i+1,3}=0\},$
we discuss the offspring distributions of $\mathcal B_{i+1,1}$ particle. Note
that an excursion $\mathcal B_{i+1,1}$ begins with a jump $\{i+1\rightarrow i-1\}$ and
ends with a jump $\{i\rightarrow i+1\}.$ So after jumping down from $i+1$ to
$i-1,$ it must have an excursion returning to $i.$ There are two
approaches for the walk to return to $i,$ that is, jumping from $i-1$
to $i$ or jumping from $i-2$ to $i$. From this point of view, one
knows that a $\mathcal B_{i+1,1}$ particle gives birth to a $\mathcal C_{i,1}$ or
$\mathcal C_{i,2}$ particle at $i$ with probability $1,$ and the approach the
walk jumping from below $i$ to $i$ determines which particle will be
born. Precisely, one has that
$$P_\omega^0(C_{i,1}+C_{i,2}=1\big|U_{i+1}=\mathbf e_4)=1.$$ Due to the same reason as (\ref{unab}),
\begin{equation}\label{unbc}\no
  P_{\omega}^0(C_{i,1}=1\big|U_{i+1}=\mathbf e_4)=1-P_{\omega}^0(C_{i,2}=1\big|U_{i+1}=\mathbf e_4)=\frac{\gamma_{i,1}}{\gamma_{i,1}+\gamma_{i,2}}.
\end{equation}

After giving birth to a $\mathcal C_{i,1}$ or $\mathcal C_{i,2}$ particle, a
$\mathcal B_{i+1,1}$ particle may give births to certain number of $\mathcal A_{i,1},$
$\mathcal A_{i,2},$ $\mathcal B_{i,1}$ and $\mathcal B_{i,2}$ excursions. Markov property
implies the independence of those born excursions. One has that
\begin{equation}\label{oca}
\begin{split}
  P_\omega^0(U_{i}=(a,b,&0,c,d,0,1,0,0)\big|U_{i+1}=(0,0,0,1,0,...,0))\\
  &=\frac{(a+b+c+d)!}{a!b!c!d!}\alpha_{i,1}^a\alpha_{i,2}^b\beta_{i,1}^c\beta_{i,2}^d(1-\alpha_{i,1}-\alpha_{i,2}-\beta_{i,1}-\beta_{i,2})\frac{\gamma_{i,1}}{\gamma_{i,1}+\gamma_{i,2}},
\end{split}
\end{equation}
 and
  \begin{equation}\label{ocb}
\begin{split}
  P_\omega^0(U_{i}=(a,b,&0,c,d,0,0,1,0)\big|U_{i+1}=(0,0,0,1,0,...,0))\\
  &=\frac{(a+b+c+d)!}{a!b!c!d!}\alpha_{i,1}^a\alpha_{i,2}^b\beta_{i,1}^c\beta_{i,2}^d(1-\alpha_{i,1}-\alpha_{i,2}-\beta_{i,1}-\beta_{i,2})\frac{\gamma_{i,2}}{\gamma_{i,1}+\gamma_{i,2}}.
\end{split}
\end{equation}
Analogously, the offspring distributions of a $\mathcal B_{i+1,3}$ particle are
\begin{equation}\label{occ}
\begin{split}
  P_\omega^0(U_{i}=(a,b,&0,c,d,0,1,0,0)\big|U_{i+1}=(0,...,0,1,0,0,0))\\
  &=\frac{(a+b+c+d)!}{a!b!c!d!}\alpha_{i,1}^a\alpha_{i,2}^b\beta_{i,1}^c\beta_{i,2}^d(1-\alpha_{i,1}-\alpha_{i,2}-\beta_{i,1}-\beta_{i,2})\frac{\gamma_{i,1}}{\gamma_{i,1}+\gamma_{i,2}},
\end{split}
\end{equation}
 and
  \begin{equation}\label{ocd}
\begin{split}
  P_\omega^0(U_{i}=(a,b,&0,c,d,0,0,1,0)\big|U_{i+1}=(0,...,0,1,0,0,0))\\
  &=\frac{(a+b+c+d)!}{a!b!c!d!}\alpha_{i,1}^a\alpha_{i,2}^b\beta_{i,1}^c\beta_{i,2}^d(1-\alpha_{i,1}-\alpha_{i,2}-\beta_{i,1}-\beta_{i,2})\frac{\gamma_{i,2}}{\gamma_{i,1}+\gamma_{i,2}}.
\end{split}
\end{equation}

{\bf \noindent(d) Offspring distribution of $\mathcal B_{i+1,2}$ particles}

\begin{center}
  \includegraphics[height=5cm]{offspringda.1}
\includegraphics[height=6cm]{offspringdb.1}
\end{center}

Next conditioned on $U_{i+1}=(0,0,0,0,1,0,0,0,0),$ that is,
$\{B_{i+1,2}=1,A_{i+1,1}=A_{i+1,2}=A_{i+1,3}=B_{i+1,1}=B_{i+1,3}
=C_{i+1,1}=C_{i+1,2}=C_{i+1,3}=0\},$ we consider the offspring
distributions of $\mathcal B_{i+1,2}$ particles.

Special attention should be payed to $\mathcal B_{i+1,2}$ excursions. An
excursion $\mathcal B_{i+1,2}$ begins with a jump $\{i+1\rightarrow i-1\}$ and ends with
a jump $\{i-1\rightarrow i+1\}.$ The point is whether it visited $i$ between
the first  and the last step. If it did not visit $i$ before the
last step, then the excursion gives birth to a $\mathcal C_{i,3}$ particle
with probability 1 and generates no any other particle. If it did
visit $i$ before the last step, that is, conditioned on
$\{C_{i,3}=0\},$ things get much more complicated. Conditioned on
$\{C_{i,3}=0\},$ after jumping down from $i+1$ to $i-1,$ the walk reaches $i$ from below $i$ at least one
time, that is, from $i-1$ to $i$ or from $i-2$ to $i.$ Therefore one
has that
$$P_\omega^0(C_{i,1}+C_{i,2}=1\big|U_{i+1}=\mathbf e_5,C_{i,3}=0)=1,$$
and similarly as (\ref{unab}) that
\begin{equation}\no
  P_\omega^0(C_{i,1}=1\big|U_{i+1}=\mathbf e_5,C_{i,3}=0)=1-P_\omega^0(C_{i,2}=1\big|U_{i+1}=\mathbf e_5,C_{i,3}=0)=\frac{\gamma_{i,1}}{\gamma_{i,1}+\gamma_{i,2}}.
\end{equation}

Since the last step of excursion $\mathcal B_{i+1,2}$ is from $i-1$ to $i+1,$ conditioned on $\{U_{i+1}=\mathbf e_5,C_{i,3}=0\},$ after reaching $i$ from below, it must jump down from $i$ to $i-1$ or from $i$ to $i-2$ at least one time, before the last step $\{i-1\rightarrow i+1\}$ occurs.
In the other words,
$$P_\omega^0(A_{i,3}+B_{i,3}=1\big|U_{i+1}=\mathbf e_5,C_{i,3}=0)=1$$ and
\begin{equation}\no
  P_\omega^0(A_{i,3}=1\big|U_{i+1}=\mathbf e_5,C_{i,3}=0)=1-P_\omega^0(B_{i,2}=1\big|U_{i+1}=\mathbf e_5,C_{i,3}=0)=\frac{\alpha_{i,3}}{\alpha_{i,3}+\beta_{i,3}}.
\end{equation}
On the other  hand, conditioned on $\{U_{i+1}=\mathbf e_5,C_{i,3}=0\},$ besides giving birthes to the above discussed particles, it may give birthes to a number of $\mathcal A_{i,1},$ $\mathcal A_{i,2},$ $\mathcal B_{i,1}$ and $\mathcal B_{i,2}$ particles.

Next we calculate the probabilities of $\{C_{i,3}=1\}$ and
$\{C_{i,3}=0\}$ conditioned on $\{U_{i+1}=\mathbf e_5\}.$

One follows from the above discussions that the sum of the products
of transition probabilities of all possible paths of an excursion
$\mathcal B_{i+1,2}$ is $\beta_{i+1,2}$ and all those possible paths could be
 divided into two classes. Paths of the first class  never
visited $i$ and paths of the second class visited $i$ from below
certain times. Sum of the product of transition probabilities of all
possible first kind paths is $\gamma_{i,3}$ and that of the second
kind paths is $\beta_{i+1,2}-\gamma_{i,3}.$ Therefore,
\begin{equation}\no
  P_\omega^0(C_{i,3}=1\big|U_{i+1}=\mathbf e_5)=1-P_\omega^0(C_{i,3}=0\big|U_{i+1}=\mathbf e_5)=\frac{\gamma_{i,3}}{\beta_{i+1,2}}.
\end{equation}

Also the Markov property implies the independence of all excursions born to $\mathcal B_{i+1,2}$ particle at $i.$ One follows by path decomposition and independence that
\begin{equation}\label{oda}P_\omega^0(U_{i}=(0,...,0,1)\big|U_{i+1}=(0,0,0,0,1,0,0,0,0))=\frac{\gamma_{i,3}}{\beta_{i+1,2}},\end{equation}
\begin{equation}\label{odb}
\begin{split}
  P_\omega^0&(U_{i}=(a,b,1,c,d,0,1,0,0)\big|U_{i+1}=(0,0,0,0,1,0,0,0,0))\\
  &=\frac{(a+b+c+d)!}{a!b!c!d!}\alpha_{i,1}^a\alpha_{i,2}^b\beta_{i,1}^c\beta_{i,2}^d(1-\alpha_{i,1}-\alpha_{i,2}-\beta_{i,1}-\beta_{i,2})\frac{(\beta_{i+1,2}-\gamma_{i,3})\gamma_{i,1}\alpha_{i,3}}{\beta_{i+1,2}(\gamma_{i,1}+\gamma_{i,2})(\alpha_{i,3}+\beta_{i,3})},
\end{split}
\end{equation}
\begin{equation}\label{odc}
\begin{split}
  P_\omega^0&(U_{i}=(a,b,1,c,d,0,0,1,0)\big|U_{i+1}=(0,0,0,0,1,0,0,0,0))\\
  &=\frac{(a+b+c+d)!}{a!b!c!d!}\alpha_{i,1}^a\alpha_{i,2}^b\beta_{i,1}^c\beta_{i,2}^d(1-\alpha_{i,1}-\alpha_{i,2}-\beta_{i,1}-\beta_{i,2})\frac{(\beta_{i+1,2}-\gamma_{i,3})\gamma_{i,2}\alpha_{i,3}}{\beta_{i+1,2}(\gamma_{i,1}+\gamma_{i,2})(\alpha_{i,3}+\beta_{i,3})},
\end{split}
\end{equation}
\begin{equation}\label{odd}
\begin{split}
  P_\omega^0&(U_{i}=(a,b,0,c,d,1,1,0,0)\big|U_{i+1}=(0,0,0,0,1,0,0,0,0))\\
  &=\frac{(a+b+c+d)!}{a!b!c!d!}\alpha_{i,1}^a\alpha_{i,2}^b\beta_{i,1}^c\beta_{i,2}^d(1-\alpha_{i,1}-\alpha_{i,2}-\beta_{i,1}-\beta_{i,2})\frac{(\beta_{i+1,2}-\gamma_{i,3})\gamma_{i,1}\beta_{i,3}}{\beta_{i+1,2}(\gamma_{i,1}+\gamma_{i,2})(\alpha_{i,3}+\beta_{i,3})},
\end{split}
\end{equation}
and
\begin{equation}\label{ode}
\begin{split}
  P_\omega^0&(U_{i}=(a,b,0,c,d,1,0,1,0)\big|U_{i+1}=(0,0,0,0,1,0,0,0,0))\\
  &=\frac{(a+b+c+d)!}{a!b!c!d!}\alpha_{i,1}^a\alpha_{i,2}^b\beta_{i,1}^c\beta_{i,2}^d(1-\alpha_{i,1}-\alpha_{i,2}-\beta_{i,1}-\beta_{i,2})\frac{(\beta_{i+1,2}-\gamma_{i,3})\gamma_{i,2}\beta_{i,3}}{\beta_{i+1,2}(\gamma_{i,1}+\gamma_{i,2})(\alpha_{i,3}+\beta_{i,3})}.
\end{split}
\end{equation}

Summing up the discussions of this section, we have the following theorem, which has also been stated as Theorem \ref{thbranc} in the introduction section.

\begin{theorem}\label{thbran}
   $\{U_{i}\}_{i\le 1}$ is a 9-type
non-homogeneous branching process with immigration  distribution as in (\ref{inia}), (\ref{inib}) and (\ref{inic}) above, and offsprings distributions as in  (\ref{oaa}-\ref{oad}), (\ref{oba}-\ref{obd}), (\ref{oca}-\ref{ocd}) and (\ref{oda}-\ref{ode}) above.
\end{theorem}

\begin{corollary}
  Let $Q_i$ be a $9\times 9$ matrix, whose $l$-th row are the means of
number of particles born  to a type-$l$ particle of the $i+1$-th
generation. The matrices $Q_i$ are called the mean matrices of the branching process $\{U_i\}_{i\le1}.$ Let
$x_i=\frac{\alpha_{i,1}}{1-\alpha_{i,1}-\alpha_{i,2}-\beta_{i,1}-\beta_{i,2}},$
$y_i=\frac{\alpha_{i,2}}{1-\alpha_{i,1}-\alpha_{i,2}-\beta_{i,1}-\beta_{i,2}},$
$z_i=\frac{\beta_{i,1}}{1-\alpha_{i,1}-\alpha_{i,2}-\beta_{i,1}-\beta_{i,2}},$
$w_i=\frac{\beta_{i,2}}{1-\alpha_{i,1}-\alpha_{i,2}-\beta_{i,1}-\beta_{i,2}},$
 $1-v=\frac{\gamma_{i,3}}{\beta_{i+1,2}},$
 $s_i=\frac{\alpha_{i,3}}{\alpha_{i,3}+\beta_{i,3}}$ and
 $t_i=\frac{\gamma_{i,1}}{\gamma_{i,1}+\gamma_{i,2}}.$ Then one calculates from the branching mechanism of $\{U_i\}_{i\le1}$ that
\begin{equation}\label{q}
  Q_i=\left(
    \begin{array}{ccccccccc}
      x_i & y_i & 0 & z_i & w_i & 0 & 0 & 0 & 0 \\
       x_i & y_i & s_i & z_i & w_i & 1-s_i & 0 & 0 & 0 \\
       x_i & y_i & 0 & z_i & w_i & 0 & 0 & 0 & 0 \\
      x_i & y_i & 0 & z_i & w_i & 0 & t_i & 1-t_i  & 0 \\
      x_iv_i & y_iv_i & s_iv_i & z_iv_i & w_iv_i & (1-s_i)v_i & t_iv_i & (1-t_i) v_i & 1-v_i \\
    x_i & y_i & 0 & z_i & w_i & 0 & t_i & 1-t_i  & 0\\
      x_i & y_i & 0 & z_i & w_i & 0 & 0 & 0 & 0 \\
       x_i & y_i & s_i & z_i & w_i & 1-s_i & 0 & 0 & 0 \\
      x_i & y_i & 0 & z_i & w_i & 0 & 0 & 0 & 0 \\
    \end{array}
  \right).
\end{equation}
\end{corollary}

\section{The  ladder time $T_1$ and the branching process--Proof of Theorem \ref{tu}}\label{prtu}
Recall that $T_1=\inf\{n\ge 0:X_n > 0\}$ is the first hitting time
of $[1,\infty).$ In this section, we aim at expressing $T_1$ in
terms of the multitype branching process $\{U_i\}_{i\le1}.$

Define \begin{equation}\no\begin{split}
  &D_{i,1}=\#\{\text{steps by the walk from above } i-1\text{ to } i-1 \text{ before time }T_1\},\\
   &V_{i,1}=\#\{\text{steps by the walk from } i-1 \text{ to } i\text{ before time }T_1\},\\
  &V_{i,2}=\#\{\text{steps by the walk from } i-2 \text{ to } i\text{ before time }T_1\}.
\end{split}
\end{equation}
Then $\sum_{i\le0}D_{i,1}$ counts all steps jumping downward by the walk before time $T_1$
and $\sum_{i\le0}V_{i,1}+V_{i,2}$ counts all steps jumping upward by the walk before time $T_1.$
But both $\sum_{i\le0}D_{i,1}$ and $\sum_{i\le0}V_{i,1}+V_{i,2}$ did not count the last step by the walk hitting $[1,\infty).$ Therefore $$T_1=1+\sum_{i\le0}D_{i,1}+V_{i,1}+V_{i,2}.$$ This together with the fact  $D_{i,1}=A_{i,1}+A_{i,2}+A_{i,3}+C_{i,1}+C_{i,2}+C_{i,3},$ $V_{i,1}=A_{i,1}+B_{i,1}+C_{i,1}$ and $V_{i,2}=A_{i,2}+B_{i,2}+C_{i,2}$ implies that
\begin{equation}\label{t}\begin{split}
  T_1&=1+\sum_{i\le0}2A_{i,1}+2A_{i,2}+A_{i,3}+B_{i,1}+B_{i,2}+2C_{i,1}+2C_{i,2}+C_{i,3}\\
  &=1+\sum_{i\le0}U_{i}(2,2,1,1,1,0,2,2,1)^T
\end{split}
\end{equation}
which proves Theorem \ref{tu}.

Next we calculate the moment of $T_1$ by mean of the branching process.

One calculates from the immigration distributions
(\ref{inia}-\ref{inic}) that
\begin{equation}\label{u1}E_\omega^0(U_1)=\z(\frac{\alpha_{1,1}}{\alpha_{1,1}+\alpha_{1,2}+\alpha_{1,3}},\frac{\alpha_{1,2}}{\alpha_{1,1}+\alpha_{1,2}+\alpha_{1,3}},\frac{\alpha_{1,3}}{\alpha_{1,1}+\alpha_{1,2}+\alpha_{1,3}},0,...,0\y)=:u_1.\end{equation}
 Therefore  one follows from Markov property that, for $i\le 0,$ $$E_\omega^0(U_i)=u_1Q_0\cdots Q_{i}.$$
 Substituting to (\ref{t}) one has that
 \begin{equation}\label{et}
 E_\omega^0(T_1)=1+\sum_{i\le0}u_1Q_0\cdots Q_i(2,2,1,1,1,0,2,2,1)^T.
 \end{equation}
\begin{remark} \label{deg}
{ \rm About the branching structure, we have the following remarks:

\begin{itemize}
  \item[-]The authors found that the branching structure  could be simplified. Indeed, note that $\mathcal{C}_{i,j}$ has the same offspring distribution with $\mathcal{A}_{i,j},$ $j=1,2,3,$ and that $\mathcal A_{i,j}$ and $\mathcal{C}_{i,j},$ $j=1,2,3,$ play the same role in the (\ref{t}). One could treat $\mathcal A_{i,j}$ and $\mathcal{C}_{i,j},$ $j=1,2,3,$ as the same type particles. In this point of view, a $6$-type branching process is enough to count exactly all steps by the walk before $T_1.$
However we still use the original $9$-type branching process because it is more understandable and each of the $9$-type particles correspond to specific jump of the walk.
\item[-] While $R=1$ the branching structure coincides with the one for (2-1) random walk constructed in Hong-Wang \cite{hw10}. Indeed, for (2-1) model, there are only three type excursion, that is $\mathcal{A}_{i,1},$  $\mathcal{B}_{i,1}$ and  $\mathcal{C}_{i,1}.$ But as the above discussed, one could treat $\mathcal{A}_{i,1}$ and $\mathcal{C}_{i,1}$ particles as the same type.
    So one need only to consider a $2$-type branching process $U_i=(A_{i,1},B_{i,1}).$
    For (2-1) random walk, one follows easily that $$\alpha_{i,1}=\omega_{i}(-1) \text{ and }\beta_{i,1}=\omega_i(-2).$$ Therefore, the distribution of $\mathcal{A}_{i+1,1}$ in (\ref{oaa}) degenerates to
    \begin{equation}\label{oaad}
\begin{split}
  P_\omega^0(U_{i}=(a,b)\big|U_{i+1}=(1,0))
  =\frac{(a+b)!}{a!b!}(\omega_i(-1))^a(\omega_{i}(-2))^b\omega_{i}(1).
\end{split}
\end{equation}

Since there is no jump of size $2$ above, one see from Figure 5 that with probability $1,$ a type $\mathcal C_{i,1}$ particle will be born to a $\mathcal B_{i+1,1}$ particle. But we treat type $\mathcal{C}_{i,1}$ particle as type $\mathcal A_{i,1}$ particle now. Therefore  the  offspring distribution of $\mathcal{B}_{i,1}$ in  (\ref{oca}) degenerates to
    \begin{equation}\label{ocad}
\begin{split}
  P_\omega^0(U_{i}=(a+1,b)\big|U_{i+1}=(0,1))
  =\frac{(a+b)!}{a!b!}(\omega_i(-1))^a(\omega_{i}(-2))^b\omega_{i}(1).
\end{split}
\end{equation}
The offspring distributions in (\ref{oaad}) and (\ref{ocad}) coincides those offspring distributions in  Hong-Wang \cite{hw10}.
\item[-] For (1-2) random walk, our result coincides with Hong-Zhang \cite{hzh10}. In this case, one needs only a $3$-type branching process, that is $U_i=(A_{i,1},A_{i,2},A_{i,3}).$
The offspring distributions of particle $\mathcal A_{i+1,1}$ and $\mathcal A_{i+1,3}$ in (\ref{oaa}) and (\ref{oab}) degenerate to
\begin{equation}\label{oaaf1}
\begin{split}
  P_\omega^0(U_{i}=(a,b,0)\big|U_{i+1}=(1,0,0))
  =\frac{(a+b)!}{a!b!}\alpha_{i,1}^a\alpha_{i,2}^b(1-\alpha_{i,1}-\alpha_{i,2}),
\end{split}
\end{equation}
\begin{equation}\label{oaaf2}
\begin{split}
  P_\omega^0(U_{i}=(a,b,0)\big|U_{i+1}=(0,0,1))
  =\frac{(a+b)!}{a!b!}\alpha_{i,1}^a\alpha_{i,2}^b(1-\alpha_{i,1}-\alpha_{i,2}).
\end{split}
\end{equation}
For (1-2) random walk, one see from Figure 4 that, except generating offsprings as $\mathcal{A}_{i+1,1}$ particle, a type $\mathcal{A}_{i,2}$ particle gives birth with probability $1$ to a type $\mathcal A_{i,3}$ particle. Therefore the offspring distribution in (\ref{oba}) degenerates to
\begin{equation}\label{obaf}
\begin{split}
  P_\omega^0(U_{i}=(a,b,1)\big|U_{i+1}=(0,1,0))
  =\frac{(a+b)!}{a!b!}\alpha_{i,1}^a\alpha_{i,2}^b(1-\alpha_{i,1}-\alpha_{i,2}).
\end{split}
\end{equation}
One see that the above (\ref{oaaf1}), (\ref{oaaf2}) and (\ref{obaf}) coincide with those offspring distributions in Hong-Zhang \cite{hzh10}.\qed
 \end{itemize}
}\end{remark}


\section{An example for testing the branching structure}\label{test}
In this section, we let $\omega_0=(q_2,q_1,p_1,p_2)$ where
$p_1,p_2,q_1,q_2>0$ and $q_2+q_1+p_1+p_2=1$ and let
$\omega=(...,\omega_0,\omega_0,\omega_0,...).$ Now consider the
random walk $\{X_n\}$ in the environment $\omega.$ Note that the
transition probabilities of $\{X_n\}$ are now independent of the
position.

We always assume that
\begin{equation}\label{pqc}p_1+2p_2-q_1-2q_2\ge0
\end{equation} which implies that $\limsup_{n\rto}X_n=\infty.$
We also mention that for such degenerated $\omega,$ $X_n$ could be
realized by i.i.d. sum. Precisely, let $S_{n}:=\sum_{i=1}^n\xi_i,$
where $\{\xi_n\}$ is an independent  sequence of random variables
with common distribution $P(\xi_1=1)=p_1,$ $P(\xi_1=2)=p_2,$
$P(\xi_1=-1)=q_1$ and $P(\xi_1=-2)=q_2.$ Then
$\{X_n\}\overset{\mathcal{D}}{=}\{S_n\}.$

 Let \begin{equation}\label{m}\no
  M=\left(
    \begin{array}{ccc}
      -\frac{q_1+q_2}{q_2} & \frac{p_1+p_2}{q_2} & \frac{p_2}{q_2} \\
     1  &  0& 0\\
      0& 1 & 0 \\
    \end{array}
  \right).
\end{equation}
 One follows from (\ref{pqc}) that $M$ has three simple
eigenvalues. Let $f,$ $g$ and $h$ be the eigenvalues of $M$ such
that $|f|>|g|>|h|.$ Then one follows also from (\ref{pqc}) that
$|f|>|g|\ge1$ and $-1<h<0.$

Indeed, define $F(\lambda):=|\lambda-M
E|=\lambda^3+\frac{q_1+q_2}{q_2}\lambda^2-\frac{p_1+p_2}{q_2}\lambda-\frac{p_2}{q_2}.$
Then $F(0)=-\frac{p_2}{q_2}<0,$ and $F(-1)=\frac{1-q_2-p_2}{q_2}>0.$
Therefore $F(x)$ has a root  $h\in (-1,0).$ Note also that $F(y)>0$
for $y$ large enough and $F(1)=\frac{q_1+2q_2-(p1+2p_2)}{q_2}\le0.$
Then $F(\lambda)$ has a root in $[1,\infty).$ Also $F(-y)<0$ for all $y$ large enough. Since $F(-1)>0,$ $F(\lambda)$ has a root in $(-\infty,-1).$

 Note that the the exit probabilities in (\ref{pb1})
degenerate to
\begin{equation}\label{cpb1}\left\{\begin{array}{l}
\mcp_{b-1}(b)=\frac{\mathbf e_1M^{b-a-1}[\mathbf e_2-\mathbf e_3]^T\Big(1+\displaystyle\sum_{l=a+1}^{b-1}\mathbf e_1M^{b-l}l\mathbf e_1^T\Big)
-\mathbf e_1M^{b-a-1}\mathbf e_1^T\displaystyle\sum_{l=a+1}^{b-1}\mathbf e_1
M^{b-l}[\mathbf e_2-\mathbf e_3]^T}{\mathbf e_1M^{b-a-1}\mathbf e_1^T\displaystyle\sum_{l=a+1}^{b-1}\mathbf e_1
M^{b-l}[\mathbf e_1-\mathbf e_2]^T-\mathbf e_1M^{b-a-1}[\mathbf e_1-\mathbf e_2]^T\Big(1+\displaystyle\sum_{l=a+1}^{b-1}\mathbf e_1M^{b-l}\mathbf e_1^T\Big)},\\
\mcp_{b-2}(b)=\frac{\mathbf e_1M^{b-a-1}[\mathbf e_2-\mathbf e_3]^T\displaystyle\sum_{l=a+1}^{b-1}\mathbf e_1M^{b-l}[\mathbf e_1-\mathbf e_2]^T
-\mathbf e_1M^{b-a-1}[\mathbf e_1-\mathbf e_2]^T\displaystyle\sum_{l=a+1}^{b-1}\mathbf e_1
M^{b-l}[\mathbf e_2-\mathbf e_3]^T}{\mathbf e_1M^{b-a-1}\mathbf e_1^T\displaystyle\sum_{l=a+1}^{b-1}\mathbf e_1
M^{b-l}[\mathbf e_1-\mathbf e_2]^T-\mathbf e_1M^{b-a-1}[\mathbf e_1-\mathbf e_2]^T\Big(1+\displaystyle\sum_{l=a+1}^{b-1}\mathbf e_1M^{b-l}\mathbf e_1^T\Big)}.\end{array}\right.
\end{equation}

Letting $a\rightarrow-\infty,$ by some careful calculations, one
follows from (\ref{cpb1}) that
\begin{equation}\label{f012}
\begin{split}
&f_{b-1}(1):=\lim_{a\rightarrow-\infty}\mcp_{b-1}({b})=1+h;\\
&f_{b-2}(1):=\lim_{a\rightarrow-\infty}\mcp_{b-2}({b})=1+h+h^2.
\end{split}
\end{equation}
We remark that, for $g=1,$ and $g\neq1,$  it is a bit different to find
the limit of (\ref{cpb1}) as $a\rightarrow-\infty.$ But the limits
have the same form for both $g=1$ and $g\neq1.$

Recall that $\mcp_{b-1}({b})$ is the simplification of
$\mcp_{b-1}(a,b,b)$ which by definition equals to
$$P_\omega^{b-1}(\text{the walk exits the interval }[a+1,b-1]\text{ at
}b).$$ Therefore one sees from ({\ref{f012}) that
\begin{equation}\begin{split}\label{f0}\no
&f_{b-1}(1)=P_\omega^{b-1}(\text{the walk exits the interval
}(-\infty,b-1]\text{ at }b)=1 +h;\\
&f_{b-2}(1)=P_\omega^{b-2}(\text{the walk exits the interval
}(-\infty,b-1]\text{ at }b)=1+h+h^2.
\end{split}\end{equation}
Similarly as all $M_k$ degenerate to $M,$ one follows by some
careful calculations from (\ref{pb2}) that
\begin{equation}\begin{split}\label{f1}\no
&f_{b-1}(2)=P_\omega^{b-1}(\text{the walk exits the interval
}(-\infty,b-1]\text{ at }b+1)=-h;\\
&f_{b-2}(2)=P_\omega^{b-2}(\text{the walk exits the interval
}(-\infty,b-1]\text{ at }b+1)=-h-h^2.
\end{split}\end{equation}

Recall that $T_1:=\inf\{n>0:X_n>0\}.$ Then
$$f_0(1)=P_\omega^0(X_{T_1}=1)=1+h$$ and
$$f_0(2)=P_\omega^0(X_{T_1}=2)=-h.$$
Therefore one has that
\begin{equation}\label{et1p}
  E_\omega^0(X_{T_1})=1f_0(1)+2f_0(2)=1-h.
\end{equation}
Since $h\in(-1,0),$ $1<E_\omega^0(X_{T_1})<2,$ which is also natural
by intuition.

One should note that the above calculations of the mean of $X_{T_1}$
involve only the exit probabilities of the walk from $(-\infty,0]$.

On the other hand, recall that $T_1$ could be expressed by the
9-type branching process constructed  in Section \ref{bran}. That
is, \begin{equation}\no
  T_1=1+\sum_{i\le0}U_{i}(2,2,1,1,1,0,2,2,1)^T.
\end{equation}

Since  $\omega_i=\omega_0$ for all $i,$ we omit the subscript $``i"$
in the notation related. For example we write $\alpha_{i,1}$ as
$\alpha_1,$ $w_i$ as $w$ et al. One has the following results:
$\alpha_1=-q_1p_1h/p_2,$
 $\alpha_{2}=q_1(p_2+hp_1+hp_2)/p_2,$
 $\alpha_3=-q_1h;$
 $\beta_1=-q_2p_1h(1+h)/p_2,$
$\beta_{2}=q_2(p_2+(p_1+p_2)h(1+h))/p_2,$
 $\beta_3=-q_2h(1+h);$
 $\gamma_1=-q_2p_1h/p_2,$
 $\gamma_{2}=q_2(p_2+hp_1+hp_2)/p_2,$
 $\gamma_3=-q_2h;$
 $x=q_1p_1h^2/p_2^2,$
 $y=-q_1h(p_2+hp_1+hp_2)/p_2^2,$
 $z=q_2p_1h^2(1+h)/p_2^2,$
 $w=-q_2h^3(1+q_2h)/p_2^2;$
 $v=1+p_2/(h(1+q_2h)),$
 $s=q_1/(q_1+q_2(1+h)),$
 $t=-p_1h/(p_2(1+h));$
$u_1=(-p_1h/p_2,(p_2+p_1h+p_2h)/p_2,-hp_2/p_2,0,0,0,0,0,0).$ With the
above notations, the mean matrix $Q_i$ of the branching process
$\{U_i\}$ in (\ref{q}) degenerates to
\begin{equation}\label{cq}
  Q=\left(
    \begin{array}{ccccccccc}
      x & y & 0 & z & w & 0 & 0 & 0 & 0 \\
       x & y & s & z & w & 1-s & 0 & 0 & 0 \\
       x & y & 0 & z & w & 0 & 0 & 0 & 0 \\
      x & y & 0 & z & w & 0 & t & 1-t  & 0 \\
      xv & yv & sv & zv & wv & (1-s)v & tv & (1-t) v & 1-v\\
    x & y & 0 & z & w & 0 & t & 1-t  & 0\\
      x & y & 0 & z & w & 0 & 0 & 0 & 0 \\
       x & y & s & z & w & 1-s & 0 & 0 & 0 \\
      x & y & 0 & z & w & 0 & 0 & 0 & 0 \\
    \end{array}
  \right).
\end{equation}

Then one has from (\ref{et}) that
 \begin{equation}\label{ccet}\no
 E_\omega^0(T_1)=1+\sum_{i\le0}u_1Q^{i+1}(2,2,1,1,1,0,2,2,1)^T.
 \end{equation}
 Note that if one assumes \begin{equation}\no
E_\omega^0(X_1)=p_1+2p_2-q_1-2q_2>0,
 \end{equation} then
 Ward Equation \begin{equation}\label{ward}
E_\omega^0(X_{T_1})=E_\omega^0(T_1)E_\omega^0(X_1)
 \end{equation}
 should hold.

Recall that $E_\omega^0(X_{T_1})$ was calculate in (\ref{et1p}) by
the exit probability and that $E_\omega^0(T_1)$ was calculated by
the branching structure constructed in Section \ref{bran}. Therefore
it will provide a good testification of the branching structure to
show that Ward Equation (\ref{ward}) holds. Equivalently, to show
(\ref{ward}), one needs only to show
\begin{equation}\label{ev}
1+\sum_{i\le0}u_1Q^{i+1}(2,2,1,1,1,0,2,2,1)^T=\frac{1-h}{p_1+2p_2-q_1-2q_2}.
\end{equation}
Some elementary calculation shows that $Q$ has four nonzero
eigenvalues $\lambda_1,$ $\lambda_2,$ $\lambda_3,$ $\lambda_4,$ and
all other eigenvalues are $0.$ One could find the eigenvectors to get a
matrix $B$ such that
\begin{equation}\label{qs}
Q=B\Lambda B^{-1}\end{equation} with $$\Lambda=\left(
                                         \begin{array}{cccccc}
                                           \lambda_1 &  &  &  & &\\
                                            &  \ddots& &  &  &  \\
                                           &  & \lambda_4&   &  &  \\
                                            & & &0 &  &  \\
                                             & & & & \ddots &  \\
                                           &  & &  & & 0 \\
                                         \end{array}
                                       \right).
$$
Theoretically, one could substitute (\ref{qs}) to the left-hand
side of (\ref{ev}) to show that (\ref{ev}) does hold. But the
calculations are technical and tedious. Instead of giving such tedious
calculations, some numerical test  may be preferred. By the
``Matlab" we get the following table.

\noindent\begin{tabular}{|c|c|c|c|c|c|c|c|}\hline
\multicolumn{8}{|c|}{\rule[-3mm]{0mm}{8mm}\bfseries Test of the
branching structure}\\ \hline $p_1$&
$p_2$&$q_1$&$q_2$&$E_\omega^0(X_1)$&$\begin{array}{c}
                                        E_\omega^0(X_{T_1})\\
:=E_\omega^0(T_1)E_\omega^0(X_1)
                                     \end{array}
$ &$\begin{array}{c}
    E_\omega^0(X_{T_1})\\
:=f_0(1)+2f_0(2)
 \end{array}
$ & Error\\
\hline 0.2100&0.3500&0.3600&0.0800&0.3900&1.467727692&1.467727692&-6.6613e-016\\
\hline 0.3000&0.2100&0.3000&0.1900&0.0400&1.323718710&1.323718710&  2.8644e-014\\
\hline 0.1789&0.3211&0.1801&0.3199&0.0012&1.481684406&1.481684406&1.8585e-013\\
\hline 0.4998&0.0002&0.4999&0.0001&0.0001&1.000399840&1.000399840&-3.5456e-011\\
\hline 0.3627&0.1373&0.3628&0.1372&0.0001&1.226498171&1.226490265&7.9058e-006\\
\hline
\end{tabular}

The column ``$E_\omega^0(X_{T_1})=E_\omega^0(T_1)E_\omega^0(X_1) $ "
of the table means that
$$E_\omega^0(X_{T_1})=E_\omega^0(T_1)E_\omega^0(X_1)=E_\omega^0(T_1)(p_1+2p_2-q_1-2q_2),$$
where $E_\omega^0(T_1)$ is calculated by the multitype branching process $\{U_i\},$ that is,
$$E_\omega^0(T_1)=1+\sum_{i\le0}u_1Q^{i+1}(2,2,1,1,1,0,2,2,1)^T.$$

\section{Invariant measure equation and law of large numbers of (2-2) RWRE--Proof of Theorem \ref{lln}}
In this section, we consider random walk $\{X_n\}$ in random
environment $\omega.$ Since by assumption $E^0(T_1)<\infty,$ then $P^0$-a.s.,
 $T_1<\infty.$

Define $\overline{\omega}(n)=\theta^{X_n}\omega.$ The process
$\{\overline{\omega}(n)\}$ is called the environment viewed from
particles. One easily show that $\{\om(n)\}$ is indeed a Markov
process under either $P_\omega^0$ or $P^0,$ with transitional kernel
\begin{equation}\label{trk}\no
K(\omega,d\omega')=\omega_0(2)\delta_{\theta^2\omega=\omega'}
+\omega_0(1)\delta_{\theta\omega=\omega'}
+\omega_0(-1)\delta_{\theta^{-1}\omega=\omega'}
+\omega_0(-2)\delta_{\theta^{-2}\omega=\omega'}.
\end{equation}
It is important to find the invariant measure and the corresponding
invariant density for the transition kernel $K(\omega,d\omega').$
If one has the invariant density in the hand, then one
could show the law of large numbers for random walk in random
environment $\{X_n\}$ and the limit velocity of the transient walk
could be expressed by the invariant density.

We borrow some notations from \cite{hzh10} to give the invariant
measure. Define
$\varphi_{\theta^k\omega}^1=P^0_{\theta^k\omega}(X_{T_1}=1)$ and
$\varphi_{\theta^k\omega}^2=P^0_{\theta^k\omega}(X_{T_1}=2).$ Whenever
$E^0(T_1)<\infty,$ define
\begin{equation}\no
  Q(d\omega)=E^0\z(\frac{1_{X_{T_1}=1}}{\varphi_{\omega}^1}\sum_{i=0}^{T_1-1}1_{\om(i)\in d\omega}
  +\frac{1_{X_{T_1}=2}}{\varphi_{\omega}^2}\sum_{i=0}^{T_1-1}1_{\om(i)\in
  d\omega}\y)\text{ and }
  \overline{Q}(d\omega)=Q(d\omega)/Q(\Omega).
\end{equation}
Then one follows verbatim as \cite{hzh10} that the measure $Q$ is
invariant under transition kernel $K(\omega,d\omega').$ Precisely,
one has that, for $B\in \mathcal{F},$
$$Q(B)=\iint1_{\omega'\in B}K(\omega,d\omega')Q(d\omega).$$
Also, following \cite{hzh10} one shows that
\begin{equation}\label{ivd}
  \frac{dQ}{dP}=\sum_{i\le0}E_{\theta^{-i}\omega}^0(N_i\big|X_{T_1}=1)+E_{\theta^{-i}\omega}^0(N_i\big|X_{T_1}=2),
\end{equation}
where $N_i:=\#\{k\in[0,T_1):X_{k}=i\}.$

The branching structure enables us to calculate the right-hand side
of (\ref{ivd}) and give specifically $\frac{dQ}{dP}.$ In fact, note
that for $i\le -2,$
\begin{equation}\no\begin{split}
  N_i=&A_{i+1,1}+A_{i+1,2}+A_{i+1,3}+C_{i+1,1}+C_{i+1,2}+C_{i+1,3}\\
  &+A_{i,1}+A_{i,2}+B_{i,1}+B_{i,2}+C_{i,1}+C_{i,2}.
\end{split}\end{equation}
Then
\begin{equation}\begin{split}\label{eni}
 E_\omega^0(N_i|X_{T_1}=1)=&\frac{E_\omega^0(U_{i+1}(1,1,1,0,0,0,1,1,1)^T;X_{T_1}=1)}{P_\omega^0(X_{T_1}=1)}\\
 &\hspace{1cm}+\frac{E_\omega^0(U_{i}(1,1,0,1,1,0,1,1,0)^T;X_{T_1}=1)}{P_\omega^0(X_{T_1}=1)}.
\end{split}
  \end{equation}
  Temporally, we set $\mathbf v_1=(1,1,1,0,0,0,1,1,1)^T$ and $\mathbf v_2=(1,1,0,1,1,0,1,1,0)^T.$
The first term in the right-hand side of (\ref{eni}) equals to
\begin{equation}\no\begin{split}
  \frac{1}{f_0(1)}&E_\omega^0(U_{i+1}\mathbf v_1;A_{1,1}=1)+\frac{1}{f_0(1)}E_\omega^0(U_{i+1}\mathbf v_1;A_{1,2}=1)\\
&=\frac{\omega_1(-1)}{\alpha_{1,1}+\alpha_{1,2}}\Big(E_\omega^0(U_{i+1}\mathbf v_1\big|A_{1,1}=1)P_\omega^0(A_{1,1}=1)+E_\omega^0(U_{i+1}\mathbf v_1\big|A_{1,2}=1)P_\omega^0(A_{1,2}=1)\Big)\\
&=\frac{\omega_1(-1)}{\alpha_{1,1}+\alpha_{1,2}}\Big(\frac{\alpha_{1,1}}{\alpha_{1,1}+\alpha_{1,2}+\alpha_{1,3}}\mathbf e_1Q_0\cdots
Q_{i+1}\mathbf v_1+\frac{\alpha_{1,2}}{\alpha_{1,1}+\alpha_{1,2}+\alpha_{1,3}}\mathbf e_2Q_0\cdots
Q_{i+1}\mathbf v_1\Big)\\
&=\Big(\frac{\alpha_{1,1}}{\alpha_{1,1}+\alpha_{1,2}},\frac{\alpha_{1,2}}{\alpha_{1,1}+\alpha_{1,2}},0,...,0\Big)Q_0\cdots
Q_{i+1}\mathbf v_1,
\end{split}
\end{equation}
since $f_0(1)=(\alpha_{1,1}+\alpha_{1,2})/\omega_1(-1)$ and
$\alpha_{1,1}+\alpha_{1,2}+\alpha_{1,3}=\omega_{1}(-1).$

Similarly the second term in the right-hand side of (\ref{eni})
equals to
$$\Big(\frac{\alpha_{1,1}}{\alpha_{1,1}+\alpha_{1,2}},\frac{\alpha_{1,2}}{\alpha_{1,1}+\alpha_{1,2}},0,...,0\Big)Q_0\cdots
Q_{i}\mathbf v_2.$$

Therefore,
 \begin{equation}\label{enif}\begin{split}
 E_\omega^0(N_i|X_{T_1}=1)=\Big(\frac{\alpha_{1,1}}{\alpha_{1,1}+\alpha_{1,2}},\frac{\alpha_{1,2}}{\alpha_{1,1}+\alpha_{1,2}},0,...,0\Big)(Q_0\cdots
Q_{i+1}\mathbf v_1+Q_0\cdots Q_{i}\mathbf v_2).
\end{split}
  \end{equation}
One the other hand,
\begin{equation}\label{enif1}\no\begin{split}
 E_\omega^0(N_i|X_{T_1}=2)= E_\omega^0(N_i|A_{1,1}=2)=\mathbf e_3(Q_0\cdots
Q_{i+1}\mathbf v_1+Q_0\cdots Q_{i}\mathbf v_2).
\end{split}
\end{equation}
Next since $N_0=1+A_{0,1}+A_{0,2}+B_{0,1}+B_{0,2},$ one has that
\begin{equation}\label{01}\begin{split}
E_\omega^0(N_0\big|X_{T_1}=2)&=E_\omega^0(N_0\big|A_{1,3}=1)\\
&=1+\mathbf e_3Q_0(1,1,0,1,1,0,0,0,0)^T=1+\mathbf e_3Q_0\mathbf v_2;\\
E_\omega^0(N_0\big|X_{T_1}=1)
&=1+\Big(\frac{\alpha_{1,1}}{\alpha_{1,1}+\alpha_{1,2}},\frac{\alpha_{1,2}}{\alpha_{1,1}+\alpha_{1,2}},0,...,0\Big)Q_0(1,1,0,1,1,0,0,0,0)^T\\
&=1+\Big(\frac{\alpha_{1,1}}{\alpha_{1,1}+\alpha_{1,2}},\frac{\alpha_{1,2}}{\alpha_{1,1}+\alpha_{1,2}},0,...,0\Big)Q_0\mathbf v_2,
\end{split}\end{equation}
  where the second equality holds due to the special structure of the mean matrix $Q_0.$
Note that
$$N_{-1}=A_{0,1}+A_{0,2}+A_{0,3}+A_{-1,1}+A_{-1,2}+B_{-1,1}+B_{-1,2}+C_{-1,1}+C_{-1,2},$$
where we mention that if the last step before $T_1$ is $\{-1\rightarrow1\},$
$A_{0,3}=1,$ otherwise, $A_{0,3}=0.$

Similarly as above, one has that
\begin{equation}\label{1}\begin{split}
 E_\omega^0(N_{-1}|X_{T_1}=1)&=\Big(\frac{\alpha_{1,1}}{\alpha_{1,1}+\alpha_{1,2}},\frac{\alpha_{1,2}}{\alpha_{1,1}+\alpha_{1,2}},0,...,0\Big)(Q_0
(1,1,1,0,...,0)^T+Q_0Q_{-1}\mathbf v_2)\\
&=\Big(\frac{\alpha_{1,1}}{\alpha_{1,1}+\alpha_{1,2}},\frac{\alpha_{1,2}}{\alpha_{1,1}+\alpha_{1,2}},0,...,0\Big)(Q_0
(1,1,1,0,0,0,1,1,1)^T+Q_0Q_{-1}\mathbf v_2);\\
 E_\omega^0(N_{-1}|X_{T_1}=2)&=E_\omega^0(N_{-1}|A_{1,3}=1)=\mathbf e_3(Q_0
(1,1,1,0,...,0)^T+Q_0Q_{-1}\mathbf v_2)\\
&=E_\omega^0(N_{-1}|A_{1,3}=1)=\mathbf e_3(Q_0
(1,1,1,0,0,0,1,1,1)^T+Q_0Q_{-1}\mathbf v_2)
\end{split}
  \end{equation}
Substituting (\ref{enif}), (\ref{01}) and (\ref{1}) to (\ref{ivd}) one concludes
that

\begin{equation}  \label{invd}
\begin{split}
\frac{dQ}{dP}=2&+\sum_{i=0}^\infty\Big(\frac{\alpha_{1+i,1}}{\alpha_{1+i,1}+\alpha_{1+i,2}},\frac{\alpha_{1+i,1}}{\alpha_{1+i,1}+\alpha_{1+i,2}},1,0,...,0\Big)Q_i\cdots
Q_0(1,1,0,1,1,0,1,1,0)^T\\
&+\sum_{i=1}^\infty\Big(\frac{\alpha_{1+i,1}}{\alpha_{1+i,1}+\alpha_{1+i,2}},\frac{\alpha_{1+i,1}}{\alpha_{1+i,1}+\alpha_{1+i,2}},1,0,...,0\Big)Q_i\cdots
Q_1(1,1,1,0,0,0,1,1,1)^T\\&=:\Pi(\omega).
  \end{split}
\end{equation}

We mention that the purpose of deriving the invariant measure
$Q(d\omega)$ and the invariant density $dQ/dP$ is to prove a law of
large number for the (2-2) random walk in random environment by an
approach known as ``the environment viewed from particles".

One follows similarly as Zeitouni \cite{ze04} Corollary 2.1.25 that
$\{\om(n)\}$ is stationary and ergodic under the measure
$\overline{Q}\otimes P_\omega^0.$ Define the local drift at site $x$
in the environment $\omega$ as $d(x, \omega)= E_{\omega}^x (X_1-x).$
The ergodicity of $\{\omega(n)\}$ under $Q \otimes P_\omega^0$ implies
that:
$$\frac{1}{n}\sum_{k=0}^{n-1}d(X_k,\omega)=\frac{1}{n}\sum_{k=0}^{n-1}d(0,\om(k))
\overset{n\rto}{\longrightarrow}E_{\overline{Q}}(d(0,\omega))\quad
\overline{Q}\otimes P_\omega^0\text{-a.s..}$$ But
\begin{equation}\no
  \begin{split}
X_n=\sum_{i=1}^{n}(X_{i}-X_{i-1})&=\sum_{i=1}^{n}(X_{i}-X_{i-1}-d(X_i,\omega))+\sum_{i=1}^nd(X_i,\omega)\\
&=:M_n+\sum_{i=1}^nd(X_i,\omega).
  \end{split}
\end{equation}
Similarly as Zeitouni \cite{ze04}, $\{M_n\}$ is a martingale and
$P^0$-a.s.,
$$\lim_{n\rto}\frac{M_n}{n}=0.$$
Therefore $$\lim_{n\rto}X_n=E_{\overline{Q}}(d(0,\omega))=:V_P.$$
Next we calculate $V_P.$ Note that
\begin{equation}\label{vp}\begin{split}
V_P&=E_{\overline{Q}}(d(0,\omega))=E_{\overline{Q}}(X_1)=E_{\mathbb P}\Big(\Pi(\omega)(2\omega_0(-2)+\omega_0(-1)+\omega_0(1)+2\omega_0(2))\Big)\Big/Q(\Omega)\\
&=\frac{E_{\mathbb P}\Big(\Pi(\omega)(2\omega_0(-2)+\omega_0(-1)+\omega_0(1)+2\omega_0(2))\Big)}{E^0(T_1\big|X_{T_1}=1)+E^0(T_1\big|X_{T_1}=2)}
\end{split}\end{equation}
But the denominator in (\ref{vp}) equals to
\begin{equation}\label{dom}E_{\mathbb P}(2+\Big(\frac{\alpha_{1,1}}{\alpha_{1,1}+\alpha_{1,2}},\frac{\alpha_{1,2}}{\alpha_{1,1}+\alpha_{1,2}},1,0,...,0\Big)\sum_{i\le0}Q_0\cdots Q_i(2,2,1,1,1,0,2,2,1)^T):=E_{\mathbb P}(D(\omega)).\end{equation}
Substituting to (\ref{vp}) one has that
$$V_P=\frac{E_{\mathbb P}\Big(\Pi(\omega)(2\omega_0(-2)+\omega_0(-1)+\omega_0(1)+2\omega_0(2))\Big)}{E_{\mathbb P}(D(\omega))}.$$

 \noindent{\large{\bf Acknowledgements:}} The authors would like to
thank Ms. Hongyan Sun and Lin Zhang for their stimulating
discussions.

{\center
\end{document}